\documentclass[11pt]{article}
\usepackage{amssymb,epsfig,amsmath}
\usepackage[dvips]{color}
\usepackage[T1]{fontenc}
\usepackage[francais]{babel}
\usepackage[utf8]{inputenc}
\usepackage{comment}
\usepackage{hyperref}
\usepackage{graphicx}
\usepackage{multicol}
\usepackage{caption}
\newtheorem{defi}{Definition}
\newtheorem{prop}[defi]{Proposition}
\newtheorem{theo}[defi]{Theorem}
\newtheorem{conj}[defi]{Conjecture}
\newtheorem{lemm}[defi]{Lemma}
\newtheorem{coro}[defi]{Corollary}
\newtheorem{rema}[defi]{Remark}
\newtheorem{exem}[defi]{Example}
\newtheorem{exems}[defi]{Examples}

\newcommand{\bdefi}{\begin{defi}}
\newcommand{\edefi}{\end{defi}}
\newcommand{\bprop}{\begin{prop}}
\newcommand{\eprop}{\end{prop}}
\newcommand{\btheo}{\begin{theo}}
\newcommand{\etheo}{\end{theo}}
\newcommand{\blemm}{\begin{lemm}}
\newcommand{\brema}{\begin{rema}}
\newcommand{\erema}{\end{rema}}
\newcommand{\bexer}{\begin{exem}}
\newcommand{\eexer}{\end{exem}}
\newcommand{\bexems}{\begin{exems}}
\newcommand{\eexems}{\end{exems}}
\newcommand{\bconj}{\begin{conj}}
\newcommand{\econj}{\end{conj}}
\newcommand{\elemm}{\end{lemm}}
\newcommand{\bcoro}{\begin{coro}}
\newcommand{\ecoro}{\end{coro}}
\newcommand{\dem}{\noindent{\bf Proof. }}
\newcommand{\rem}{\noindent{\bf Remark. }}

\usepackage{mathrsfs}
\renewcommand\mathcal{\mathscr}

\newcommand{\T}{{\cal T}}

\newcommand{\M}{{\cal M}}

\newcommand{\G}{{\cal G}}

\newcommand{\C}{{\cal C}}

\newcommand{\maths}[1]{{\mathbb #1}}  

\newcommand{\RR}{\maths{R}}
\newcommand{\NN}{\maths{N}}

\newcommand{\SSS}{\maths{S}}

\newcommand{\ZZ}{\maths{Z}}

\newcommand{\cqfd}{\hfill$\Box$}

\newcounter{fig}

\newcommand{\dddp}{\partial_{\infty}^2\widetilde{\Sigma}}
\newcommand{\ddp}{\partial_{\infty}\widetilde{\Sigma}}

\newcommand{\Image}{\operatorname{Image}}

\newcommand{\Supp}{\operatorname{Supp}}

\newcommand{\CAT}{\operatorname{CAT}}

\newcommand{\revet}{(\widetilde{\Sigma},[\widetilde{q}])}

\newcommand{\srfce}{({\Sigma},{[q]})}

\newcommand{\revetm}{(\widetilde{\Sigma},\widetilde{m})}
\newcommand{\srfcem}{({\Sigma},m)}
\newcommand{\Id}{\operatorname{Id}}

\newcommand{\Lq}{\Lambda_{[q]}}

\newcommand{\Lqr}{\widetilde{\Lambda}_{[\widetilde{q}]}}

\newcommand{\Lmr}{\widetilde{\Lambda}_{\widetilde{m}}}
\newcommand{\Gqr}{\G_{[\widetilde{q}]}}
\newcommand{\Gmr}{\G_{\widetilde{m}}}
\newcommand{\Gq}{\G_{[q]}}

\newcommand{\Ltilde}{\widetilde{\Lambda}}
\newcommand{\mutilde}{\widetilde{\mu}}
\newcommand{\grperevet}{\Gamma_{\widetilde{\Sigma}}}

\newcommand{\BP}{\operatorname{FS}}

\title{Measured flat geodesic laminations.}
\author{Thomas~Morzadec}

\begin{document}
\maketitle

Département de mathématique, UMR 8628 CNRS, Université Paris-Sud, Bât. 430, F-91405 Orsay Cedex, France  
Bureau : 16.  {\it thomas.morzadec@math.u-psud.fr}

\medskip

\textbf{Abstract: } Since their introduction by Thurston, measured geodesic laminations on hyperbolic surfaces occur in many contexts. 
In \cite{Morzy1}, we have introduced a notion of flat laminations on surfaces endowed with a half-translation structure (that is a singular flat surface with holonomy 
$\{\pm\Id\}$), similar to
geodesic laminations on hyperbolic surfaces. 
Here is a sequel to this article that aims at defining transverse measures on flat laminations similar to
transverse measures on hyperbolic laminations, taking into account that two different leaves of a flat lamination may no longer be disjoint. 
One aim of this paper is to construct a tool that could allow a fine description of the space of degenerations of half-translation structures on a surface. In this paper, 
we define a nicer topology than the Hausdorff topology on the set of measured flat laminations and a natural continuous 
projection of the space of measured 
flat laminations onto the space of
measured hyperbolic laminations, for some arbitrary half-translation structure and hyperbolic metric on a surface. We prove in particular that the space of measured flat
laminations is projectively compact.   
\footnote{ Keywords: Measured geodesic lamination, flat surface, half-translation structure, holomorphic quadratic differential, 
measured foliation, hyperbolic surface, dual tree. 
AMS codes 30F30, 53C12, 53C22.}

%

\section{Introduction.}

The main aim of this article is to propose a definition of transverse measure on the (geodesic) flat laminations, introduced in \cite{Morzy1}, on a surface
endowed with a half-transla\-tion structure, that is a flat metric with conical singular points and with holonomies in $\{\pm\Id\}$. Although the definition is inspired of
transverse measures on the geodesic laminations on hyperbolic surfaces (see for instance \cite{Bonahon97}), the extension is non trivial, notably since the images of
two leaves of a 
flat lamination are not necessarly disjoint. We will call {\it measured flat lamination} a flat lamination endowed with a transverse measure. We will define a sufficiently
fine
topology on the set of measured flat laminations and  we will construct a (non injective) natural continuous projection of the space of measured flat laminations onto the
space of
measured
hyperbolic laminations, for any choice of a half-translation structure and of a (complete) hyperbolic metric on a surface, and we will describe its lack of injectivity.
This allows to consider the 
measured flat laminations that are the limits of some sequences of periodic local geodesics, in the projectivized space of measured flat laminations. 
This in turn could yield a better understanding of the degenerations of half-translation structures on a surface, as initiated in \cite{DucLeiRaf10}. 
In particular, as spaces of measures are suitable for analysis tools (distributions as in \cite{Bonahon97}), this could allow a finer study of the boundary of the space of
half-translation structures that we will develop in a subsequent work. We refer to \cite{Morzy3}
for a survey of \cite{Morzy1} and of this work.

We use the same notations as in \cite{Morzy1}: let $\Sigma$ be a compact, connected, orientable surface, without boundary (to simplify in the introduction).
A {\it half-translation structure} (or flat structure with conical singularities and holonomies in $\{\pm\Id\}$) on $\Sigma$ is the data consisting in a 
(possibly empty) discrete set of points $Z$ of $\Sigma$ and of a
Euclidean metric on $\Sigma-Z$ with conical singular points of angles of the form $k\pi$, with $k\in\NN$ and $k\geqslant 3$ at each point of $Z$, such that the holonomy of every
piecewise $\C^1$ loop of $\Sigma-Z$ is contained in $\{\pm\Id\}$. We refer to Section \ref{structuresplates} notably when the boundary is not empty.

The surface $\Sigma$ endowed with a half-translation structure is a complete and locally $\CAT(0)$ metric space $(\Sigma,d)$. Let $p:(\widetilde{\Sigma},\widetilde{d})\to
(\Sigma,d)$ be a locally isometric universal cover. Two local geodesics ${\ell},{\ell}'$ of $(\Sigma,d)$, defined up to 
changing the
origins, are said to be 
{\it interlaced} if they have some lifts $\widetilde{\ell},\widetilde{\ell}'$ in $\widetilde{\Sigma}$ such that the image of $\widetilde{\ell}$ intersects both
complementary
components of $\widetilde{\ell}'(\RR)$ in $\widetilde{\Sigma}$, and conversely. A local geodesic is said to be {\it self-interlaced} if it is interlaced with itself. 
We endow the set of oriented, but non parametrized, local geodesics of $(\Sigma,d)$ with the quotient topology of the compact-open topology for the action by translations
on the parametrizations,
of $\RR$ on the parametrized local geodesics, that is called the {\it geodesic topology}. 

\bdefi\cite[Déf.~2.2]{Morzy1} A (geodesic) flat lamination on $(\Sigma,d)$ is a non empty set $\Lambda$ of complete local geodesics of $(\Sigma,d)$, defined up to changing 
origin, whose elements are called {\it leaves}, such that:

\begin{itemize}
 \item[$\bullet$]the leaves of $\Lambda$ are non self-interlaced and pairwise non interlaced;
 \item[$\bullet$] $\Lambda$ is invariant by changing the orientations of the leaves;
 \item[$\bullet$] $\Lambda$ is closed for the geodesic topology.
\end{itemize}

We will call {\it support} of $\Lambda$ the union of the images of the leaves of $\Lambda$. 
\edefi

New phenomenons appear in flat laminations compared with hyperbolic ones: the images of two leaves are generally not disjoint, the flat laminations are
not determined by their supports (uncountably many flat laminations can have the same support), the cylinder components (see Section \ref{geodesiclaminations}) may have uncontably many leaves.
Finally, there are three types of minimal components of a flat lamination on a compact surface (periodic leaf travelled in both orientations, minimal component of 
recurrent type or of finite graph type, see Theorem \ref{1Morzy} for a complete statement). Compared with hyperbolic laminations, the main difficulty to define
transverse 
measures on flat laminations is that the images of the leaves are not necessarly disjoint and that the support does not determine the lamination. Hence, we no longer define
the transverse measure
as a family of measures on the images of the arcs transverse to the lamination, but as a family of measures on the sets of local geodesics that  intersect them transversally,
and we have to refine the notion of invariance by holonomy of these families of measures.

\medskip

In the first section, we define the measured flat laminations and we endow their set with a topology. In the second one, we define the preimage of a measured
flat lamination in
a cover. In the third one, we define a homeomorphism between the space of measured flat laminations on a compact surface endowed with a fixed half-translation structure
and the space of Radon measures on the set of geodesics (defined up to changing origin) of a locally isometric universal cover, that are invariant by the covering
group action and whose supports are some flat laminations. In the fourth one, we define  a proper, surjective, continuous map from the space of measured flat laminations
to the space of measured hyperbolic laminations, for a fixed complete hyperbolic metric with totally geodesic boundary, and we characterize its lack of injectivity. In the fifth
one, we define the intersection number between a measured flat lamination and a free homotopy class of closed curves. In the last one, we define the  tree associated to a 
measured flat lamination together with the universal covering group action on it. These tools should be useful to do analysis on the space of degenerations of half-translation
structures on surfaces, and we plan to develop this in future work (see \cite{Morzy4}).

\medskip\noindent
{\small{\it Acknoledgement: I want to thank Frederic Paulin for many advices and corrections that have deeply improved the redaction of this paper.}

\section{Definitions.}\label{rappelsconventions}

In this section, we recall the definition and some properties of half-translation structures on surfaces and the definition of flat laminations introduced in \cite{Morzy1}.
Then, we define measured flat laminations and we endow their set with a topology.

\subsection{Half-translation structures.}\label{structuresplates}

As in \cite{Morzy1}, in the whole paper, we will use the definitions and notation of \cite{BriHae99} for a surface endowed with a distance
$(\Sigma,d)$: (locally) $\CAT(0)$, $\delta$-hyperbolic,... 
Notably, a {\it geodesic} (resp. a {\it local geodesic}) 
of $(\Sigma,d)$ is an isometric (resp. locally isometric) map $\ell:I\to \Sigma$, where $I$ is an interval of $\RR$. It will be called a 
  {\it segment}, a {\it ray} or a {\it  geodesic line} of $(\Sigma,d)$ if $I$ is respectively a compact interval, a closed half line 
  (generally $[0,+\infty[$)
  or $\RR$.
  If there is no precision, a {\it geodesic} is a geodesic line. A
 {\it germ of geodesic ray}, or simply a {\it germ}, is an equivalence class of locally geodesic rays for the equivalence relation  
 $r_1\sim_0 r_2$ if $r_1$ and $r_2$ coïncide on a non empty initial segment that is not reduced to a point. Similarly, the relation $r\sim_\infty r'$ if there exist $T,T'>0$
 such that $r(t+T)=r'(t+T')$ for all $t\geqslant 0$, is an equivalence relation on the set of subrays of a local geodesic.
An equivalence class for this equivalence relation is called an {\it end} (in the sense of Freudhental) of a local geodesic. A local geodesic has two ends. We call  
{\it geodesic topology} the compact-open topology on the set $\G_d$ of local geodesics for the distance $d$ 
   or the quotient topology of this topology by the action by translations of $\RR$ at the source, on the set $[\G_d]$ of local geodesics defined
   up to changing origin. The quotient map from $\G_d$ to $[\G_d]$  will be
   denoted by $g\mapsto[g]$, and if $F$ is a subset of $\G_d$, we will denote by $[F]$ its image in $[\G_d]$.
    
    \medskip

Let $\Sigma$ be a connected, orientable surface, with (possibly empty) boundary. 
Assume that $\Sigma$ is endowed with a Euclidean metric on $\Sigma-Z$, where $Z$ is a discrete subset of $\Sigma$.
If the holonomy of every piecewise $\C^1$ loop in $\Sigma-Z$ is contained in $\{\pm\Id\}$, two vectors $v_1$ and $v_2$ tangent to 
$\Sigma$ are said to have the 
{\it same direction}  if $v_2$ is the image of $\pm v_1$ by holonomy along a piecewise $\C^1$ path in $\Sigma-Z$ between the basepoints of 
$v_1$ and $v_2$. 
This definition does not depend on the choice of this path. A piecewise $\C^1$ path or union of paths is said to have  {\it constant direction}, if all its tangent vectors,
at the points in $\Sigma-Z$, have the same direction.

\bdefi
A half-translation structure (or flat structure with conical singularities and holonomies in $\{\pm\Id\}$) on a surface $\Sigma$ is the data of a (possibly empty) discrete
subset $Z$ of $\Sigma$
and a Euclidean metric on $\Sigma-Z$ with conical singularity of angle $k_z\pi$ at each $z\in Z$, with $k_z\in\NN$, $k_z\geqslant 3$
  if $z\in Z-Z\cap\partial\Sigma$ and $k_z\geqslant 2$  if $z\in Z\cap\partial\Sigma$, such that the holonomy of every piecewise $\C^1$ loop in 
 $\Sigma-Z$ is contained in $\{\pm\Id\}$ and such that the union of the boundary components has constant direction.
\edefi

 We will denote by $[q]$ a half-translation structure on $\Sigma$, with $q$ a holomorphic quadratic differential on $\Sigma$ (see \cite[§~2.5]{Morzy1} for an explanation
 of the notation
 and to \cite[Def.~1.2~p.~2]{Strebel84} for a definition of a holomorphic quadratic differential in the case where the boundary is non empty).
 A half-translation structure defines a geodesic distance $d$ on $\Sigma$ that is locally $\CAT(0)$. We will call {\it local flat geodesics} the local geodesics of a half-translation 
 structure. A continuous map $\ell:\RR\to\Sigma$ is a local flat geodesic if and only if it satisfies (see \cite[Th.~5.4~p.24]{Strebel84} and
 \cite[Th.~8.1~p.~35]{Strebel84}): for every $t\in\RR$,

\medskip
\noindent
$\bullet$~  if $\ell(t)$ does not belong to $Z$, there exists a neighborhood $V$ of $t$ in $\RR$ such that $\ell_{|V}$ is a Euclidean segment
(hence, $\ell_{|V}$ has constant direction);
 
\noindent
$\bullet$~  if $\ell(t)$ belongs to $Z-Z\cap\partial\Sigma$, then the two angles defined by the germs of $\ell([t,t+\varepsilon[)$ and $\ell(]t-\varepsilon,t])$, 
with $\varepsilon>0$ small enough, measured in both connected components of $U-\ell(]t-\varepsilon,t+\varepsilon[)$, with $U$  a small enough neighborhood of
$\ell(t)$, are at least $\pi$.

\noindent
$\bullet$~  if $\ell(t)$ belongs to $Z\cap\partial\Sigma$, then the angle defined by the germs of $\ell([t,t+\varepsilon[)$ and $\ell(]t-\varepsilon,t])$, 
with $\varepsilon>0$ small enough, measured in the connected component of $U-\ell(]t-\varepsilon,t+\varepsilon[)$ which is disjoint from $\partial\Sigma$, with $U$  a small 
enough neighborhood of $\ell(t)$, is at least $\pi$.
\begin{center}
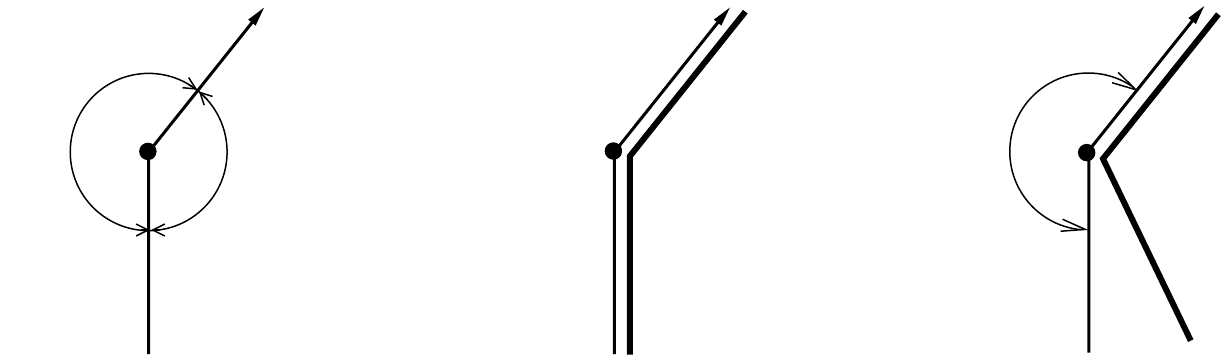
\end{center}

\subsection{Geodesic laminations on surfaces endowed with a half-translation structure and with a (complete) hyperbolic metric.}\label{geodesiclaminations}

Let $\Sigma$ be a connected, orientable surface with (possibly empty) boundary. Let $[q]$ be a half-translation structure and let $m$ be  a hyperbolic metric with totally 
geodesic boundary on $\Sigma$. Let  
$p:\widetilde{\Sigma}\to\Sigma$ be a  universal cover of covering group $\Gamma_{\widetilde{\Sigma}}$, let $[\widetilde{q}]$ be the unique half-translation structure 
and let $\widetilde{m}$ be the unique hyperbolic metric  
on $\widetilde{\Sigma}$ such that $p:\revet\to\srfce$ and $p:\revetm\to\srfcem$ are locally isometric. To be consistent with \cite{Morzy1}, and notably to be allowed to use
\cite[Rem.~2.9]{Morzy1}, we will always assume that $\Sigma$ is a cover (possibly trivial) of a compact surface whose Euler characteristic is negative.
In particular, if $\widetilde{d}$ is the distance defined by $[\widetilde{q}]$ or $\widetilde{m}$, according to the theorem of Cartan-Hadamard, the metric space 
$(\widetilde{\Sigma},\widetilde{d})$ is complete, $\CAT(0)$, and $\delta$-hyperbolic, with $\delta>0$. 
Furthermore, there exists a unique $\Gamma_{\widetilde{\Sigma}}$-equivariant
homeomorphism between the boundaries at infinity of $\widetilde{\Sigma}$ for the two
metrics, thank to which we identify them. Let $\partial_\infty\widetilde{\Sigma}$ denote this boundary at infinity and
$\dddp=\ddp\times\ddp-\{(x,x),x\in\ddp\}$.

 In \cite[§~2.3]{Morzy1}, we have given a very global definition of interlaced local geodesics, in a locally $\CAT(0)$, complete, connected metric space, 
 whose boundary at infinity of a universal cover is endowed with a (total) cyclic order. Here, we only recall the specific definition in the case of
 connected, orientable surfaces, endowed with a complete locally $\CAT(0)$ metric.  
  Since $(\widetilde{\Sigma},\widetilde{d})$ is  $\CAT(0)$, the intersection of the images of two geodesics of $(\widetilde{\Sigma},\widetilde{d})$ is connected, possibly 
  empty. If $\widetilde{\ell}$ is a geodesic of $(\widetilde{\Sigma},\widetilde{d})$, since $\widetilde{\ell}(\RR)$ is not necessarly disjoint from $\partial\widetilde{\Sigma}$, the complementary 
  $\widetilde{\Sigma}-\widetilde{\ell}(\RR)$ may have more than two connected components. However, the orientation of $\widetilde{\Sigma}$ allows 
  do distinguish the two sides 
  of $\widetilde{\ell}(\RR)$ (except if $\widetilde{\ell}(\RR)$ is a boundary component of $\widetilde{\Sigma}$, in which case $\widetilde{\ell}(\RR)$ has only one side).
  Two geodesics $\widetilde{\ell},\widetilde{\ell'}$ of $(\widetilde{\Sigma},\widetilde{d})$, defined up to changing the origins, are
{\it interlaced} if neither $\widetilde{\ell}(\RR)$ nor $\widetilde{\ell}'(\RR)$ is a boundary component of $\widetilde{\Sigma}$, and
$\widetilde{\ell}$ intersects two connected components of $\widetilde{\Sigma}-\widetilde{\ell}'(\RR)$ corresponding to one and the other of its sides
(or the same after exchanging 
$\widetilde{\ell}$ and $\widetilde{\ell}'$, which is equivalent). Two local geodesics of  
$(\Sigma,d)$ are {\it interlaced} if they admit some lifts in $(\widetilde{\Sigma},\widetilde{d})$ which are interlaced, and a local geodesic is 
{\it self-interlaced} if it is interlaced with itself. We easily convinced that this definition is equivalent to \cite[§~2.3]{Morzy1} in the case of surfaces. 

If $d$ is the distance defined by $m$, two local geodesics are non interlaced if and only if they are disjoint and a local geodesic is non self-interlaced if and only if 
it is simple. We refer to \cite[§~3.1]{Morzy1} for a characterization of the local geodesics for $[q]$ that are non interlaced.

\bdefi\label{laminationplate1}
A geodesic lamination (or simply a lamination) of $(\Sigma,d)$, with $d$ the distance defined by $m$ or $[q]$, is a non-empty set  
$\Lambda$ of (complete) local geodesics of $(\Sigma,d)$, defined up to changing origin, whose elements are called leaves, such that: 

\medskip

\noindent
$\bullet$~ leaves are non self-interlaced;

\noindent
$\bullet$~ leaves are pairwise non interlaced;

\noindent
$\bullet$~ if $\ell$ belongs to $\Lambda$ then so does $\ell^-$, with $\ell^{-}(t)=\ell(-t)$;

\noindent
$\bullet$~ $\Lambda$ is closed for the geodesic topology.
\edefi

We say that $\Lambda$ is a {\it  flat lamination} if $d$ is defined by $[q]$ and that $\Lambda$ is a {\it hyperbolic lamination} if $d$ is defined by $m$.
Usually, a hyperbolic lamination of $\srfcem$ is defined as
a non-empty closed subset of $\Sigma$, that is a union  of images of simple and pairwise disjoint local geodesics of $\srfcem$.
The definitions are equivalent in the case of hyperbolic laminations but not in the case of flat laminations 
(see \cite[§~4.1]{Morzy1}).
We recall the two main results of \cite{Morzy1} about flat laminations. In Theorem \ref{1Morzy}, a {\it cylinder component} is a maximal set of leaves of $\Lambda$ whose images are contained in a non degenerated flat cylinder (hence, 
these leaves are periodic),
a {\it minimal component} is a sublamination which is the closure, for the geodesic topology, of a leaf $\ell$ and its opposite $\ell^-$. The minimal component is of 
{\it recurrent type} if $\ell$ is regular (i.e. does not meet any singular point) and is not periodic, all the images of its leaves are then dense in a {\it domain} of $\Sigma$, i.e. the closure 
of a connected open subset bounded by some periodic local geodesics, if it is not equal to $\Sigma$. The minimal component is of {\it finite graph type} if the image of
$\ell$ is a finite graph, and if neither $\ell$ nor its 
opposite are eventually periodic. All the images of its leaves are then equal, and no leaf is eventually periodic. 
An end of a local geodesic {\it terminates} in a minimal component or in a cylinder component if there exists a ray in the equivalence class of this end which is the ray of
a leaf of the minimal component or a ray of a boundary component of the corresponding flat cylinder.

\btheo\label{1Morzy}\cite[§~6]{Morzy1} Let $\Lambda$ be a flat lamination on a compact, connected, orientable surface endowed with a half-translation structure.
 Then $\Lambda$ is a finite union of cylinder components, of minimal components (of recurrent type, finite graph type and periodic leaf travelled in both senses)
and of isolated leaves (for the geodesic topology) both of whose ends terminate in a minimal component or a cylinder component.
\etheo

A {\it cyclic orientation} on a finite metric graph $X$ is the data of cyclic orders (see \cite[§~2.3.1]{Wolf11} for the definition) on the sets of germs of 
locally geodesic rays issued from each vertice of $X$.

\btheo Every cyclically oriented, connected, finite, metric graph $X$, without extremal point, may be the support of an uncountable minimal flat lamination
with no eventually periodic leaf,  on a 
compact and connected surface endowed with a  half-translation structure, except if $X$ is homeomorphic to a circle, a dumbbell pair, a flat height or a flat theta, 
by a homeomorphism preserving the cyclic
orientation (i.e \begin{picture}(0,0)%
\includegraphics{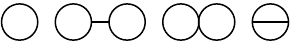}%
\end{picture}%
\setlength{\unitlength}{3771sp}%
\begingroup\makeatletter\ifx\SetFigFont\undefined%
\gdef\SetFigFont#1#2#3#4#5{%
  \reset@font\fontsize{#1}{#2pt}%
  \fontfamily{#3}\fontseries{#4}\fontshape{#5}%
  \selectfont}%
\fi\endgroup%
\begin{picture}(1461,198)(-8,830)
\end{picture}%
, where the orientations are given by the plan).
\etheo    

\subsection{Definition of measured flat laminations.}

 Let $\srfce$ be a connected, orientable surface with (possibly empty) boundary, endowed with a half-translation structure. An {\it arc} is a piecewise $\C^1$ map
 $\alpha:[0,1]\to\Sigma$ which is a homeomorphism onto its image. Let $\Lambda$ be a flat lamination of $\srfce$. An arc $\alpha$ is {\it transverse}
 to a leaf or to a segment of leaf $\ell$ of $\Lambda$ if

\noindent
$\bullet$~$\alpha$ is transverse to $\ell$ outside the singular points of $[q]$ and the singular points of $\alpha$;

\medskip
\noindent
$\bullet$~ for every singular point $x$ of $[q]$ or of $\alpha$ in $\Image(\ell)\cap\alpha(]0,1[)-\Image(\ell)\cap\alpha(]0,1[)\cap\partial\Sigma$, there exists a neighborhood
$U$ of $x$ that is a topological disk, and a segment $S$ of $\ell$ such that $U-\Image(S)\cap U$ has two connected components and the connected components of
$U\cap(\alpha([0,1])-\{x\})$ are contained in different components of $U-\Image(S)\cap U$;

\noindent
$\bullet$~ for every singular point $x$ of $[q]$ or of $\alpha$ in $\Image(\ell)\cap\alpha(]0,1[)\cap\partial\Sigma$, there exists a neighborhood
$U$ of $x$ that is a topological quarter of disk, and a segment $S$ of $\ell$ such that $U-\Image(S)\cap U$ has two or three connected components and the connected
components of
$U\cap(\alpha([0,1])-\{x\})$ are contained in two different components of $U-\Image(S)\cap U$, with only one which is not disjoint from $\partial\Sigma$. The arc $\alpha$ 
may intersect $\partial\Sigma$ along a segment;
\begin{center}
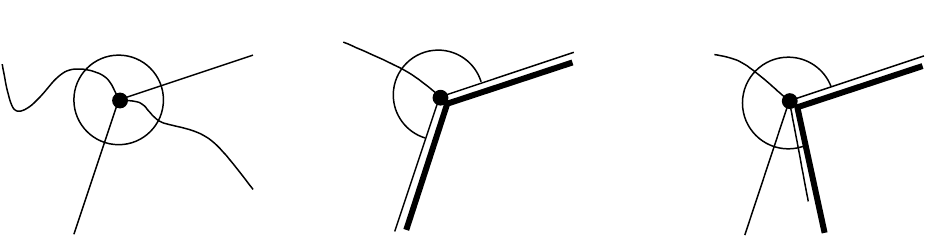
\end{center}

\noindent
$\bullet$~ $\alpha$ is tangent to $\ell$  neither in $0$  nor in $1$. However, $\alpha(0)$ and $\alpha(1)$ may belong to $\ell(\RR)$.

\medskip

An arc $\alpha$ is {\it transverse}  to a set $F$ of leaves or of segments of leaves of  $\Lambda$ if it is transverse to every element of $F$, and $F$ is  
{\it transverse to $\alpha$} if $\alpha$ is transverse to $F$. In particular, an arc is  transverse to $\Lambda$ if it is transverse to every leaf of $\Lambda$.

If $\alpha:[0,1]\to\Sigma$ is an arc of  
 $\Sigma$, we denote by $G(\alpha)$ the subset of $\Gq$ consisting in the local geodesics of $\srfce$ which are transverse to $\alpha$ and whose origins 
 belong to 
 $\alpha([0,1])$. By definition, if $\alpha'([0,1])\subseteq\alpha([0,1])$, then $G(\alpha')\subseteq G(\alpha)$. Let $F_1\subseteq\Gq$ be such that $[F_1]\subseteq\Lambda$,
 and let $\alpha_1$ and $\alpha_2$ be two disjoint  arcs transverse to $F_1$, such that
$F_1\subseteq G(\alpha_1)$ and every element of $F_1$ intersects $\alpha_2([0,1])$ at a positive time. For every $g_1\in F_1$, we define $t_{g_1}=
\min\{t>0\;:\;g_1(t)\in \alpha_2([0,1])\}$. Let $F_2$ be the subset of the elements $g_2\in G(\alpha_2)$ such that there exists $g_1\in F_1$ with 
$g_2(t)=g_1(t+t_{g_1})$ for all $t\in\RR$. A {\it holonomy} $h:F_1\to F_2$ of $\Lambda$ is a homeomorphism between $F_1$ and $F_2$ defined by 
$h(g_1)=g_2:t\mapsto g_1(t+t_{g_1})$ such that there exists a homotopy $H:[0,1]\times[0,1]\to\Sigma$ between $\alpha_1=H(\cdot,0)$ and $\alpha_2=H(\cdot,1)$ such that:

\noindent
$\bullet$~ for every $t\in\,[0,1]$,  the map $s\mapsto H(s,t)$ is an arc transverse to every segment of leaf $g_{1|[0,\,t_{g_1}]}$, with $g_1\in
F_1$;

\noindent
$\bullet$~ for every $\ell\in F_1$, there exists $s_\ell\in[0,1]$ such that $t\mapsto H(s_\ell,t)$ is a segment of $\ell$ (up to changing the parametrization);

\noindent
$\bullet$~ the intersections $H([0,1]\times]0,1[)\cap \alpha_i([0,1])$ with $i=1,2$ are empty.

Contrarily to the case of measured foliations, if the images of the geodesics are not pairwise disjoint, the map $H$ may not be injective.

\bdefi\label{defmesuretransverse} 
A transverse measure on $\Lambda$ is a family $\mu=(\mu_\alpha)_\alpha$ of Radon measures $\mu_\alpha$ defined on $G(\alpha)$, for every  arc $\alpha$ transverse to $\Lambda$, such that:

\noindent
$(1)$~ the support of $\mu_\alpha$ is the set $\{\ell\in G(\alpha)\;:\;[\ell]\in\Lambda\}$;

\noindent
$(2)$~ if $h:F_1\to F_2$ is a holonomy of $\Lambda$, where $\alpha_1,\alpha_2$ are two disjoint  arcs transverse to $F_1$ and $F_1\subset G(\alpha_1)$ and 
$F_2\subset G(\alpha_2)$ are some Borel sets, then $h_*(\mu_{\alpha_1|F_1})=\mu_{\alpha_2|F_2}$;

\noindent
$(3)$~ $\mu_\alpha$ is $\iota$-invariant, with $\iota(\ell)=\ell^-:t\mapsto\ell(-t)$;

\noindent
$(4)$~ if $\alpha'([0,1])\subseteq\alpha([0,1])$, then $\mu_{\alpha|G(\alpha')}=\mu_{\alpha'}$. 
\edefi

We will denote by $(\Lambda,\mu)$ a flat lamination endowed with a transverse measure, that we will call a {\it measured flat lamination}, and we will denote by 
$\mathcal{M}\mathcal{L}_p(\Sigma)$ the set of measured flat laminations on $\Sigma$.
We endow $\mathcal{M}\mathcal{L}_p(\Sigma)$ with the topology such that a sequence $(\Lambda_n,\mu_n)_{n\in\NN}$ converges to $(\Lambda,\mu)$ if and only if for every arc 
$\alpha$, if $\alpha$ is transverse to $\Lambda$, then $\alpha$ is transverse to $\Lambda_n$ for $n$ large enough and $\mu_{n,\alpha}\overset{*}{\rightharpoonup}\mu_\alpha$ 
in the space of Radon measures on $G(\alpha)$.

\medskip

A leaf $\ell$ of $\Lambda$ is {\it positively recurrent} if there exists an arc $\alpha$ transverse to ${\ell}$ such that $\ell$ intersects $\alpha([0,1])$ at an infinite 
number of positive times. For example, if $\Sigma$ is compact, the leaves terminating in some minimal components are positively recurrent.

\blemm\label{quedescomposnatesminimales}
If $\Lambda$ is endowed with a transverse measure $\mu$, then the only leaves of $\Lambda$ which are isolated and positively recurrent are the periodic leaves. 
\elemm

\dem Assume there exists an isolated leaf $\ell$ of $\Lambda$ which is positively recurrent, and let $\alpha$ be an arc as above.
      Then, there exists an increasing sequence $(t_n)_{n\in\NN}$ of successive positive times such that $\ell(t_n)
 \in\alpha([0,1])$.   
 For every $n\in\NN$, let $\ell_n:t\mapsto\ell(t+t_n)$. Let $n\in\NN$, and let $t_n=s_0\leqslant\dots\leqslant s_k=t_{n+1}$ be the finite sequence of times between $t_n$ and 
 $t_{n+1}$ such that $\ell(s_i)$ is a singular point for every $i\in[1,k-1]\cap\NN$, and choose some arcs $\alpha_i$ and $\beta_i$
 transverse to $\ell$ at $\ell(s_i)$ and at $\ell(s_{i+1})$ if $i\leqslant k-1$. We always can choose $\alpha_i$ and $\beta_i$ such that there exists a homotopy between 
 $\alpha_i$ and $\beta_i$ that allows to define a holonomy 
 between $\ell_{\alpha_i}:t\mapsto\ell(t+s_i)\in G(\alpha_i)$ and $\ell_{\beta_i}:t\mapsto\ell(t+s_{i+1})\in G(\beta_i)$, and such that 
 the arcs $\alpha_0$ and $\beta_{k-1}$ are some subarcs of $\alpha$
 (If $\ell(s_i)$ or $\ell(s_{i+1})$ belong to $\partial\Sigma$, in order to have a holonomy between $\alpha_i$
 and $\beta_i$, we may choose $\alpha_i$ or $\beta_i$ such that they intersect $\partial\Sigma$ along a segment or such one of their extremal points belong to 
 $\ell(\RR)$).
 
 \begin{center}
  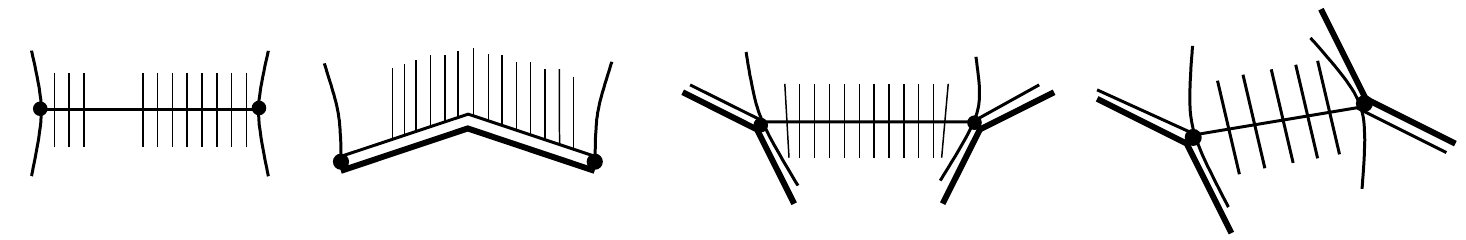
 \end{center}

 Then, according to the assertion $(2)$ of Definition \ref{defmesuretransverse} 
 we have $\mu_{\alpha_i}(\ell_{\alpha_i})=\mu_{\beta_i}(\ell_{\beta_i})$. We can assume, without loss of generality, that for every $i\in[0,k-1]\cap\NN$, 
 the intersection $\alpha_{i+1}([0,1])\cap\beta_i([0,1])$ is the image of an arc transverse to $\ell$, thus according to the assertion $(4)$ of Definition \ref{defmesuretransverse},
 we have $\mu_{\alpha_{i+1}}(\ell_{\alpha_{i+1}})=\mu_{\beta_i}(\ell_{\beta_i})$, thus $\mu_{\alpha_{i+1}}(\ell_{\alpha_{i+1}})=\mu_{\alpha_i}(\ell_{\alpha_i})$,
 and by iteration $\mu_\alpha(\ell_n)=\mu_{\alpha_0}(\ell_{\alpha_0})=\mu_{\beta_{k-1}}(\ell_{\beta_{k-1}})=\mu_\alpha(\ell_{n+1})$ (the holonomies have to be built 
 piece by piece instead of directly between $\ell_n$ and $\ell_{n+1}$, to include the case where $\ell(\RR)$ is not disjoint from $\partial\Sigma$, 
 and thus we cannot define
 the homotopy globally).
 
 By iteration, we have  $\mu_\alpha(\{\ell_n\})=\mu_\alpha(\{\ell_0\})$ for all $n\in\NN$. Since $\ell$ is isolated,   
 the leaf $\ell_0$ is isolated in $G(\alpha)$, and since $\ell_0$ belongs to the support of $\mu_\alpha$, we have $\mu_\alpha(\ell_0)>0$. Finally,  
 the set $\{\ell_n\}_{n\in\NN}$ is contained in
 $G(\alpha)$, which is relatively compact according to the theorem of Ascoli. And if $\ell$ were not periodic, the set $\{\ell_n\}_{n\in\NN}$ would be infinite and 
 $\mu_\alpha(\overline{\{\ell_n\}_{n\geqslant 0}})
 \geqslant\sum_{n\geqslant 0}\mu_\alpha(\ell_0)$ would be infinite. Hence $\mu_\alpha$ would not be locally finite, which is a contradiction.\cqfd
 
 \section{Preimage of a measured flat lamination.}\label{relevez}
 
 Let $\srfce$ be a connected, orientable surface, with (possibly empty) boundary, endowed with a half-translation structure and let $(\Lambda,\mu)$ be a measured 
 flat lamination on $\Sigma$. Let $p':(\Sigma',[q'])\to\srfce$ be a locally isometric cover of $\srfce$ with covering group  $\Gamma_{\Sigma'}$ and let 
 $\Lambda'$ be the preimage of $\Lambda$ in $\Sigma'$ (see the remark before \cite[Lem.~2.3]{Morzy1}). Since $p'$ is a local diffeomorphism, 
 if $\alpha$ is an arc of $\Sigma$ which is transverse to $\Lambda$, and if $\alpha'$ is a lift of $\alpha$ in $\Sigma'$, then $\alpha'$ is transverse to $\Lambda'$,
 and $p'$ induces a homeomorphism $f_{\alpha'}:G(\alpha')\to G(\alpha)$. We set $\mu_{\alpha'}=(f_{\alpha'}^{-1})_*\mu_\alpha$. More generally, if $\alpha'$ is
 any arc 
 transverse to $\Lambda'$, then its image is a union  of images of some lifts of some  arcs transverse to $\Lambda$, say 
 $\alpha'([0,1])=\alpha_1'([0,1])\cup\dots\cup\alpha_n'([0,1])$ such that, for every $2\leqslant k\leqslant n$, the intersection $\alpha_{k-1}'([0,1])\cap
 \alpha_{k}'([0,1])$ is the image of a lift of an arc transverse to  $\Lambda$, and if $p\not\in\{k,k+1,k-1\}$, then $\alpha_k'([0,1])\cap\alpha_{p}'([0,1])=\emptyset$. 
 Then, for every $k\in[1,n-1]\cap\NN$, we have 
 
 \begin{align*}
  \mu_{\alpha_k'|G(\alpha_k')\cap G(\alpha_{k+1}')}&=(f_{\alpha'_k}^{-1})_*(\mu_{\alpha_k|G(p'\circ\alpha_k')\cap
 G(p'\circ\alpha_{k+1}')})\\
&=(f_{\alpha'_{k+1}}^{-1})_*(\mu_{\alpha_{k+1}|G(p'\circ\alpha_k')\cap G(p'\circ\alpha_{k+1}')})\\
&=\mu_{\alpha_{k+1}'|G(\alpha_k')\cap G(\alpha_{k+1}')}
 \end{align*}
 
 Hence, there exists a unique measure $\mu_{\alpha'}$ on $G(\alpha')$ such that $\mu_{\alpha'|G(\alpha_k')}=\mu_{\alpha_k'}$ for every $k\in[1,n]\cap\NN$.
 We define $\mu'=(\mu'_{\alpha'})_{\alpha'\in\T'}$, where $\T'$ is the set of  arcs transverse to $\Lambda'$, as the family of measures just constructed.
 By naturality, $\Lambda'$ is $\Gamma_{\Sigma'}$-invariant and $\mu'$ is invariant by the action by homeomorphisms of 
 $\Gamma_{\Sigma'}$ defined by ${\gamma(\mu'_{\alpha'})_{\alpha'\in\T'}}=(\gamma_*\mu'_{\gamma^{-1}\alpha'})_{\alpha'\in\T'}$, 
 for every $\gamma\in\Gamma_{\Sigma'}$. 
 
 \blemm\label{mu'}
 The family $\mu'$ is the unique transverse measure on $\Lambda'$ such that if $\alpha'$ is the lift of an arc 
 $\alpha$ transverse to $\Lambda$, then $\mu_{\alpha'}=(f_{\alpha'}^{-1})_*\mu_\alpha$. Moreover, the map $\mu\mapsto\mu'$ from $\mathcal{M}\mathcal{L}_p(\Sigma)$ to
 $\mathcal{M}\mathcal{L}_p(\Sigma')$ thus defined is a homeomorphism between $\mathcal{M}\mathcal{L}_p(\Sigma)$ and the space of measured flat laminations of $\Sigma'$
 which are $\Gamma_{\Sigma'}$-invariant. 
 \elemm

 \dem The properties $(1)$, $(3)$ and $(4)$ of Definition \ref{defmesuretransverse} are clearly satisfied by $\mu'$, and if $h':F_1'\subset 
 G(\alpha_1')\to F_2'
 \subset G(\alpha_2')$ is a holonomy of $\Lambda'$ which is the lift of a holonomy $h:F_1\subset G(\alpha_1)\to F_2\subset G(\alpha_2)$ of $\Lambda$,    
 with $\alpha_1,\alpha_2$ two disjoint arcs transverse to $\Lambda$, then 
 $\mu_{\alpha_2'|F_2'}=(f_{\alpha_2'}^{-1})_*(\mu_{\alpha_2|F_2})=
 (f^{-1}_{\alpha_2'})_*h_*(\mu_{\alpha_1|F_1})=h'_*(f^{-1}_{\alpha_1'})_*(\mu_{\alpha_1|F_1})=h'_*(\mu_{\alpha'_1|F'_1})$. Otherwise, we see $h'$ as a concatenation
 of compositions of holonomies that are lifts of some holonomies of $\Lambda$, and we show similarly that  
 $\mu_{\alpha_2'|F_2'}=h'_*(\mu_{\alpha_1'|F_1'})$. Hence $\mu'$ is a transverse measure on $\Lambda'$. The map from 
 $\mathcal{M}\mathcal{L}_p(\Sigma)$ to $\mathcal{M}\mathcal{L}_p(\Sigma')$ defined in this way is injective by construction.
 If $(\Lambda',\mu')$ is a measured flat lamination of $(\Sigma',[q'])$ which is $\Gamma_{\Sigma'}$-invariant, then the set of projections of the leaves of $\Lambda'$
 by $p:\Sigma'\to\Sigma$ is a flat lamination $\Lambda$ of $\srfce$, and if $\alpha$ is an arc  transverse to $\Lambda$ and $\alpha'$ is a lift of $\alpha$ 
 in $\Sigma'$, then $(f_{\alpha'})_*\mu'_{\alpha'}$ is a measure on $G(\alpha)$, that does not depend on the choice of the lift, by $\Gamma_{\Sigma'}$-invariance.
 The family of measures defined in this way  satisfies clearly the properties $(1)$, $(3)$ and $(4)$ of Definition \ref{defmesuretransverse}.
 Moreover if $h:F_1\to F_2$ is a holonomy between two Borel sets of $G(\alpha_1)$ and $G(\alpha_2)$, with $\alpha_1,\alpha_2$ two disjoint arcs transverse to $\Lambda$,
 then it lifts to a holonomy between two Borel sets of $G(\alpha_1')$ and $G(\alpha_2')$, with $\alpha_1',\alpha_2'$ two lifts of $\alpha_1,\alpha_2$,
 and we have $\mu_{\alpha_2|F_2}=h_*(\mu_{\alpha_1|F_1})$. Hence $(\Lambda,\mu)$ is a measured flat lamination, and its image by the previous map is 
 $(\Lambda',\mu')$.

Hence, this map is a bijection between $\mathcal{M}\mathcal{L}_p(\Sigma)$ and the set of measured flat laminations of $(\Sigma',[q'])$ that are $\Gamma_{\Sigma'}$-invariant.
Finally, if $(\Lambda_n,\mu_n)_{n\in\NN}$ is a sequence of $\mathcal{M}\mathcal{L}_p(\Sigma)$ that converges to $(\Lambda,\mu)$, and if $(\Lambda'_n,\mu'_n)_{n\in\NN}$ and
$(\Lambda',\mu')$ are their images in $\mathcal{M}\mathcal{L}_p(\Sigma')$, then if $\alpha'$ is an arc  transverse to $\Lambda'$ which is the lift of an arc $\alpha$ transverse
to $\Lambda$, 
the arc $\alpha$ is 
transverse to $\Lambda_n$ for $n$ large enough, and since $p$ is a local diffeomorphism, $\alpha'$ is transverse to $\Lambda'_n$ for  $n$ large enough. Moreover, we have 
$\mu'_{n,\alpha'}=(f_{\alpha'}^{-1})_*\mu_{n,\alpha}\overset{*}{\rightharpoonup}(f_{\alpha'}^{-1})_*\mu_{\alpha}=\mu'_\alpha$. If $\alpha'$ is an arc transverse
to $\Lambda'$ which is not the lift of an arc  transverse to $\Lambda$, then by finite decomposition, we prove similarly that if $n$ is large enough, $\Lambda'_n$ is transverse to 
$\alpha'$ and $\mu'_{n,\alpha'}\overset{*}{\rightharpoonup}\mu'_{\alpha'}$. Hence, the map $(\Lambda,\mu)\mapsto(\Lambda',\mu')$ is continuous.
The same property holds for its inverse, hence it is a homeomorphism.\cqfd

\section{Measured flat laminations and Radon measures on the space of flat geodesics of the universal cover.}

Let $\srfce$ be a connected, orientable surface with (possibly empty) boundary, endowed with a half-translation structure and 
let $p:\revet\to\srfce$ be a locally isometric universal
cover, with covering groupe $\grperevet$.
 The local geodesics of $\revet$ are geodesics, and if $\Ltilde$ is a flat lamination of $\revet$ and if $\alpha$ is an arc  transverse to $\Ltilde$ such that every 
 geodesic
 of $G(\alpha)$ intersects $\alpha([0,1])$ only at its origin, then the map $g_\alpha:G(\alpha)\to[G(\alpha)]$ defined by $g_\alpha(g)=[g]$ is a homeomorphism.

We denote by $\M_{\grperevet}([\Gqr])$ the space of Radon measures on the space $[\Gqr]$ endowed with the geodesic topology, 
which are $\grperevet$ and $\iota$-invariant (where $\iota(\ell)=\ell^-:t\mapsto\ell(-t)$) and whose supports are $\grperevet$-invariant flat laminations,
endowed with the weak-star topology.
Let $\nu$ be an element of $\M_{\grperevet}([\Gqr])$ whose support is $\Ltilde$ and let $\alpha$ be an  arc transverse to $\Ltilde$ such that the geodesics of $G(\alpha)$ 
intersect $\alpha([0,1])$ only at their origins.  Then $\mutilde_\alpha=(g_\alpha^{-1})_*(\nu_{|[G(\alpha)]})$ is a Radon measure on $G(\alpha)$ whose support is 
$g_\alpha^{-1}(\widetilde{\Lambda}\cap[G(\alpha)])$. If $\alpha$ is an arc  transverse to $\Ltilde$, but if some geodesic of $G(\alpha)$ possibly intersects 
$\alpha([0,1])$ at several points, then we define the measure $\mutilde_\alpha$ by a finite gluing process as in Section \ref{relevez}, which is always possible
since the geodesics of $\Gqr$ are proper. Let $\mutilde_\nu$ be the family of transverse measures defined in this way. 

\bprop\label{nu} The family of measures $\mutilde_\nu$ is a transverse measure on $\Ltilde$, and the map $\nu\mapsto\mutilde_\nu$ is a  
homeomorphism between $\M_{\grperevet}([\Gqr])$ and the space of measured flat laminations on $\revet$ which are $\grperevet$-invariant. 
\eprop

\dem Unless the opposite is specified, until the end of the proof, if $\alpha$ is an arc  transverse to a flat lamination $\widetilde{\Lambda}$, we assume that the geodesics 
of $G(\alpha)$ only intersect $\alpha([0,1])$ at their origins. We may always assume this  up to shortening $\alpha$, since the geodesics of $G(\alpha)$ are proper and 
transverse to $\alpha$. 

The properties $(1)$, $(3)$, and $(4)$ of Definition \ref{defmesuretransverse} are clearly satisfied by $\mutilde_\nu$. 
If $h:F_1\to F_2$ is a holonomy of $\Ltilde$ between two Borel sets $F_1$ of
$G(\alpha_1)$ and $F_2$ of $G(\alpha_2)$, then, by the definition of the holonomies,  the sets $[F_1]$ and $[F_2]$ are equal. Hence 
 $\mutilde_{\nu,\alpha_2|F_2}=(g_{\alpha_2}^{-1})_*(\nu_{|[F_2]})=h_*(g_{\alpha_1}^{-1})_{*}(\nu_{|[F_1]})=h_*(\mutilde_{\nu,\alpha_1|F_1})$ and $\mutilde_\nu$ 
 is invariant by holonomy, so it is a transverse measure on $\Ltilde$, which is $\grperevet$-invariant by naturality. 

Assume that two measures $\nu_1$ and $\nu_2$ of $\M_{\grperevet}([\Gqr])$ define the same measured flat lamination $(\widetilde{\Lambda},\mutilde)$
by this construction. Then, if $U$ is a relatively compact open set of $[\Gqr]$, the set $U$ is contained in a finite union of open sets $(U_i)_{1\leqslant i\leqslant n}$ 
such that for every $U_i$ there exists an  arc $\alpha_i$ transverse to  $\Ltilde$, such that $U_i\subset [G(\alpha_i)]$.  We have 
$\nu_1(U_i)=\nu_2(U_i)$ for all $1\leqslant i\leqslant n$, and then $\nu_1(U)=\nu_2(U)$. Since the relatively compact open sets span the Borel tribute, 
we have $\nu_1=\nu_2$. Hence, the map $\nu\mapsto\mutilde_\nu$ is injective. 

Let us now prove it is surjective. Let $(\Ltilde,\mutilde)$ be a measured flat lamination and let ${\cal R}$ be the sets of subsets $A\subset[\G_{[\widetilde{q}]}]$ 
such that there exists a finite family of  arcs $\alpha_1,\dots,\alpha_n$, with $n\in\NN$, transverse to $\Ltilde$ (and such that the elements of $G(\alpha_i)$ intersect 
$\alpha_i([0,1])$ only at their origins), for all $1\leqslant i\leqslant n$, and some Borel subsets $A_i\subset G(\alpha_i)$ such that $A=\amalg[A_i]$. 
Then ${\cal R}$ is non empty and $\emptyset\in{\cal R}$. If $A,B$ are
two elements of ${\cal R}$, let $(\alpha_i)_{1\leqslant i\leqslant n}$ and $(\beta_j)_{1\leqslant j\leqslant k}$ be some finite
families of arcs transverse to $\Ltilde$ and
$(A_i\subseteq G(\alpha_i))_{1\leqslant i\leqslant n}$ and $(B_j\subseteq G(\beta_j))_{1\leqslant j\leqslant k}$ be some families of
Borel subsets of $\G_{[\widetilde{q}]}$
such that $A=\amalg[A_i]$ and $B=\amalg[B_j]$. For all $1\leqslant i\leqslant n$ and $1\leqslant j\leqslant k$, we replace $A_i$ by 
$A_i-A_i\cap g_{\alpha_i}^{-1}([B_j]\cap[G(\alpha_i)])$. Then $A\cup B=(\amalg[A_i])\amalg(\amalg[B_j])$ and $A-B=\amalg [A_i]$. Hence, 
the set ${\cal R}$ is a ring. Let $A\in{\cal R}$ and let $(A_i\subseteq G(\alpha_i))_{1\leqslant i\leqslant n}$, be as above. We define $\nu(A)=\overset{n}{\underset{i=1}
{\sum}}\mu_{\alpha_i}(A_i)$. Then $\nu(A)\geqslant 0$ for all $A\in{\cal R}$ and $\nu(\emptyset)=0$. Let $(A_n)_{n\in\NN}$ be a sequence of pairwise disjoint elements of
${\cal R}$ such that
$\underset{n\in\NN}\amalg A_n$ belongs to ${\cal R}$. Let $\{\alpha_i\}_{1\leqslant i\leqslant k}$ be a family of arcs transverse to $\Ltilde$ associated
to $\underset{n\in\NN}\amalg A_n$ as above. For every $n\in\NN$, there exists a family $\{A'_{n,i}\subset G(\alpha_i)\}_{1\leqslant i\leqslant k}$ of Borel subsets 
such that 
$A_n=\overset{k}{\underset{i=1}\amalg}[A'_{n,i}]$. If $n\not=m$, we have $A'_{n,i}\cap A'_{m,i}=\emptyset$. 
Hence \begin{align*}
       \nu(A)&=\overset{k}{\underset{i=1}\sum}\,\mutilde_{\alpha_i}(\underset{n\in\NN}\amalg A'_{n,i})\\
             &=\underset{n\in\NN}\sum\,\overset{k}{\underset{i=1}\sum}\,\mutilde_{\alpha_i}(A'_{n,i})\\
             &=\underset{n\in\NN}\sum\nu(A_n)
      \end{align*}
Hence, the constructed map is $\sigma$-additive on ${\cal R}$, and since $\nu(\emptyset)=0$ it
is a premeasure on ${\cal R}$. For every arc $\alpha$ transverse to $\widetilde{\Lambda}$, 
the measure $\mu_\alpha$ is locally finite, hence if $A\in{\cal R}$ has a compact closure, then 
$\nu(A)<+\infty$. Moreover, according to the theorem of Ascoli, the set of geodesics whose origins belong to a compact subset of 
$\widetilde{\Sigma}$ is compact, and since $\widetilde{\Sigma}$ is a countable union of compact sets, the space $[\G_{[\widetilde{q}]}]$ is a countable union of compact
sets, hence $\nu$ is $\sigma$-finite. According to the theorem of extension of Caratheodory, the premeasure $\nu$ can be extended to a unique measure, still denoted by $\nu$,
on the 
$\sigma$-algebra generated by ${\cal R}$, which is equal to the Borel tribute of $[\G_{[\widetilde{q}]}]$. 
According to the properties $(1)$ and $(3)$ of Definition \ref{defmesuretransverse}, 
the support of the measure $\nu$ is $\Ltilde$ 
and $\nu$ is $\iota$ and $\grperevet$-invariant by naturality. By construction, the measured flat lamination associated to $\nu$ by the previous map
is $(\Ltilde,\mutilde)$. Hence, the map $\nu\mapsto\mutilde_\nu$ is bijective.

Finally, if $(\nu_n)_{n\in\NN}$ is a sequence in $\M_{\grperevet}([\Gqr])$ which converges to $\nu$, and if $(\Ltilde_n,\mutilde_n)_{n\in\NN}$ and $(\Ltilde,\mutilde)$
are their images in $\mathcal{M}\mathcal{L}_p\,(\widetilde{\Sigma})$, then if $\alpha$ is an arc  transverse to $\Ltilde$ (such that the geodesics of $G(\alpha)$ intersect $\alpha([0,1])$
only at their origins), for $n$ large enough, the lamination 
$\Ltilde_n$ is transverse to $\alpha$
and $\mutilde_{n,\alpha}=(g_\alpha^{-1})_*(\nu_{n|[G(\alpha)]})\overset{*}{\rightharpoonup}(g_\alpha^{-1})_*(\nu_{|[G(\alpha)]})=\mutilde_\alpha$. 
 Hence, the map $\nu\mapsto\mutilde_\nu$ is continuous. Similarly, if $(\Ltilde_n,\mutilde_n)_{n\in\NN}$ is a sequence of measured flat laminations of $\widetilde{\Sigma}$
 which are $\grperevet$-invariant, that converges to  
$(\Ltilde,\mutilde)$, and if $(\nu_n)_{n\in\NN}$ and $\nu$ are their preimages in $\M_{\grperevet}([\Gqr])$, then if $f$ is a continuous map 
from $[\Gqr]$ to $\RR$ whose support is compact, there exists a finite family $\{\alpha_1,\dots,\alpha_p\}$ of arcs transverse to  
$\Ltilde$ and to $\Ltilde_n$ for $n$ large enough, such that $\Supp(f)\subset\bigcup_{1\leqslant i\leqslant p}[G(\alpha_i)]$, and for all $1\leqslant i\leqslant p$
and $n$ large enough, we have $\nu_n(f_{|[G(\alpha_i)]})=
\mutilde_{n,\alpha_i}(f\circ g_{\alpha_i|G(\alpha_i)})\longrightarrow\mutilde_{\alpha_i}(f\circ g_{\alpha_i|G(\alpha_i)})=\nu(f_{|[G(\alpha_i)]})$. Hence 
$\nu_n\overset{*}{\rightharpoonup}\nu$ and the inverse map is continuous. Hence, the map $\nu\mapsto\mutilde_\nu$ is a homeomorphism.\cqfd

\medskip

Let $\phi:\mathcal{M}\mathcal{L}_p(\Sigma)\to\M_{\grperevet}([\Gqr])$ be the composition of the map
$\mu\mapsto\mu'$ defined at Lemma \ref{mu'}, with $\Sigma'=\widetilde{\Sigma}$, with the inverse of the map
$\nu\mapsto\mutilde_\nu$ defined at Lemma \ref{nu}.

\bcoro\label{correspondance} 
 The map $\phi$ is a homeomorphism.\cqfd
\ecoro

\section{Links between measured flat laminations and measured hyperbolic laminations.}\label{liensplatshyperboliques}

 In this Section \ref{liensplatshyperboliques}, we still denote by $\srfce$ a connected, orientable surface with (possibly empty) boundary,
 endowed with a half-translation structure, and by
 $p:\revet\to\srfce$ a locally isometric universal cover whose covering group is $\Gamma_{\widetilde{\Sigma}}$. We assume that $\Sigma$ is compact and  that
 $\chi(\Sigma)<0$, and we denote by $m$ a hyperbolic metric with totally geodesic boundary on $\Sigma$ and by $\widetilde{m}$ the unique hyperbolic metric on 
 $\widetilde{\Sigma}$ such that $p:\revetm\to\srfcem$ is locally isometric. For every geodesic
 $g$ of $[\Gqr]$ or $[\Gmr]$, we denote by
 $E(g)\in\dddp$ its ordered pair of points at infinity, and if $F$ is a set of geodesics, then $E(F)=\{E(g)\;:\;g\in F\}$. 
 The space $\mathcal{M}\mathcal{L}_h(\Sigma)$ of measured hyperbolic laminations on $\srfcem$, endowed with the topology defined in \cite[p.~19]{Bonahon97},
 is homeomorphic to the space of Radon measures on $[\Gmr]$ that are 
 $\Gamma_{\widetilde{\Sigma}}$ et $\iota$-invariant, whose supports are hyperbolic laminations, denoted by $\M_{\grperevet}([\Gmr])$
(see
 \cite[Prop.~17~p.~154]{Bonahon88}). Here, we use this fact and the homeomorphism between $\M_{\grperevet}([\Gqr])$
 and $\mathcal{M}\mathcal{L}_p(\Sigma)$ defined in Corollary \ref{correspondance} to investigate the links between
 $\mathcal{M}\mathcal{L}_p(\Sigma)$ and $\mathcal{M}\mathcal{L}_h(\Sigma)$.  
 We denote by  $\varphi:[\Gqr]\to[\Gmr]$ the map associating to the geodesic $\widetilde{\ell}\in[\Gqr]$ the geodesic $\varphi(\widetilde{\ell})\in[\Gmr]$ that 
 corresponds to $\widetilde{\ell}$, i.e. such that $E(\varphi(\widetilde{\ell}))=E(\widetilde{\ell})$
  (see \cite[§4.2]{Morzy1}). Then, $\varphi$ is surjective and continuous, and a closed subset $F$ of $[\Gqr]$ is a flat lamination if and only if $\varphi(F)$ is a hyperbolic
  lamination. Moreover  $\varphi$ is proper. Indeed, if $K$ is a compact set of $[\Gmr]$ and if $(\widetilde{\ell}_n)_{n\in\NN}$ is a sequence of $\varphi^{-1}(K)$,
  then by definition of $\varphi$, up to taking a subsequence, the sequence $(E(\widetilde{\ell}_n))_{n\in\NN}$ converges in $\dddp$. Since $\revet$ is  
  $\delta$-hyperbolic (as $\Sigma$ is compact and $\chi(\Sigma)<0$, see \cite[Rem.~2.10]{Morzy1}),
  according to \cite[Lem.~2.6]{Morzy1}, up to taking a subsequence, the sequence $(\widetilde{\ell}_n)_{n\in\NN}$ converges to a geodesic $\widetilde{\ell}$ such
  that $E(\widetilde{\ell})\in E(K)$. Hence  $\widetilde{\ell}$ belongs to $\varphi^{-1}(K)$, therefore $\varphi^{-1}(K)$ is compact. 
  
  However, the map $\varphi$ is not injective. By definition, two different geodesics $\widetilde{\ell},\widetilde{\ell}'$ of $\revet$ have the same image by
  $\varphi$ if and only if they have the same ordered pair of points at infinity, and according to the flat strip theorem (see for example 
  \cite[Th.~2.13~P.~182]{BriHae99}), their images are parallel and contained in a maximal flat strip of $\revet$. Then, their projections in $\srfce$  are freely homotopic 
  periodic local geodesics, and hence their images are contained in a maximal flat cylinder (see for example \cite[Th.~2.(c)]{MS85}). The points at infinity of
  $\widetilde{\ell}$, $\widetilde{\ell}'$ and of their images  
   $\varphi(\widetilde{\ell})=\varphi(\widetilde{\ell}')$  are therefore the attractive and repulsive fixed points of an element of the covering group, 
  and the projection of $\varphi(\widetilde{\ell})$ is a closed local geodesic of $\srfcem$. Hence, the restriction of $\varphi$ to the set of geodesics of  
  $\revet$ whose projections are not periodic is an injective map whose image is the set of geodesics of $\revetm$ whose projections are not closed geodesics.
  Moreover, the preimage of a geodesic $\widetilde{\lambda}$ of $\revetm$ whose projection  is a closed local geodesic is the set of geodesics of $\revet$ having the same
  ordered pair of points at infinity than $\widetilde{\lambda}$, which may be a unique geodesic, or may be a set of parallel geodesics whose images are contained in a flat
  strip, and whose projections are periodic.      
  
  \blemm\label{surjectionpropre}
  The map $\varphi$ defines a continuous, surjective and proper map $\varphi_*$ from the space of Radon measures on $[\Gqr]$ to the space of Radon measures 
  on $[\Gmr]$. Moreover, $\varphi_*\nu$ belongs to $\M_{\grperevet}([\Gmr])$ if and only if $\nu$ belongs to $\M_{\grperevet}([\Gqr])$ and the restriction of $\varphi_*$ to
  $\M_{\grperevet}([\Gqr])$ is a surjective map $\varphi_*:\M_{\grperevet}([\Gqr])\to\M_{\grperevet}([\Gmr])$. 
  \elemm
  
 \dem The map $\varphi$ is continuous and proper, hence it defines a continuous map $\varphi_*$ from the space of Radon measures on $[\Gqr]$ to the space of
 Radon measures on $[\Gmr]$. 
  Let us first prove that $\varphi_*$ is surjective. The map $s:[\Gmr]\to[\Gqr]$ which associates, to a hyperbolic geodesic, the flat geodesic to which it corresponds
  (if it is unique) and the "middle" geodesic of the set of geodesics (contained in a flat strip)  to which it corresponds otherwise,
  is a measurable section (that is not continuous) of $\varphi$.
  Since $\varphi$ is continuous, the preimage of a compact set by $s$ is relatively compact. Hence $s$ defines a map $s_*$ from the set of Radon measures on
  $[\Gmr]$ to the space of Radon measures on $[\Gqr]$, and $\varphi_*\circ s_*=\Id$. Therefore $\varphi_*$ is surjective.
  
  Let us now prove that $\varphi_*$ is proper. The space $[\Gqr]$ is countable at infinity. Thus, there exists a sequence $(K_n)_{n\in\NN}$ of  
  compact sets such that, for all $n\in\NN$, $K_n$ is contained in the interior of $K_{n+1}$ and $\bigcup_{n\in\NN}K_n=[\Gqr]$. If $C$ is a  compact set of the space 
  of Radon 
  measures on $[\Gmr]$
  and if $K$ is a compact set of $[\Gqr]$, then the set $\{\nu(K),\nu\in(\varphi_*)^{-1}(C)\}$ is bounded by the maximum of $\{\nu(\varphi(K)),\nu\in C\}$, which is 
  finite since $C$ is compact.
  Hence, for all $n\in\NN$, the set $\{\nu_{|K_n},\nu\in(\varphi_*)^{-1}(C)\}$ is compact and if $(\nu_k)_{k\in\NN}$ is a sequence of $(\varphi_*)^{-1}(C)$, proceeding 
  by diagonal extraction, we show that there exists a subsequence, 
  still denoted by 
  $(\nu_k)_{k\in\NN}$, and a Radon measure $\nu$ on $[\Gqr]$ such that for all $n\in\NN$,
  $\nu_{k|K_n}\overset{*}{\rightharpoonup}\nu_{ |K_n}$. And, by the choice of the sequence $(K_n)_{n\in\NN}$, for all $f\in\C_c([\Gqr])$, there exists $n\in\NN$ 
  such that 
  $\Supp(f)\subset K_n$, and then
  $(\nu_k(f))_{k\in\NN}=(\nu_{k|K_n}(f))_{k\in\NN}$ converges to $\nu_{|K_n}(f)=\nu(f)$. Hence, we have $\nu_k\overset{*}{\rightharpoonup}\nu$ and $\varphi_*$  
  is proper on the space of Radon measures.

  Finally, by definition of $\varphi$, the measure $\varphi_*\nu$ belongs to $\M_{\grperevet}([\Gmr])$ if and only if $\nu$ belongs to $\M_{\grperevet}([\Gqr])$. 
  Moreover, the space $\M_{\grperevet}([\Gmr])$ is closed (see \cite[Prop.~3~et~17]{Bonahon88}), and since $\varphi_*$ is continuous, 
  its preimage is closed. Hence, the restriction of $\varphi_*$ to these spaces is continuous, proper and surjective.\cqfd
  
  \medskip
  
  The map $\varphi_*$ is not injective. Considering the lack of injectivity of $\varphi$, the preimage by $\varphi_*$ of a Radon measure on $[\Gmr]$ whose support contains no 
  geodesic whose projection in $\srfcem$ is a closed geodesic, consists in a unique Radon measure on $[\Gqr]$ whose support contains no geodesic whose projection in $\srfce$
  is a periodic local geodesic. However, the preimage by 
  $\varphi_*$ of the $\grperevet$-orbit of a Dirac measure, whose support is the $\grperevet$-orbit of a geodesic $\widetilde{\lambda}$, of mass $\delta$, whose projection
  in $\srfcem$ 
  is a closed geodesic, is the set of Radon measures on $[\Gqr]$, whose support is the $\grperevet$-orbit of a closed subset $F$
  of the set of geodesics of $\revet$ having the same ordered pair of points at infinity that $\widetilde{\lambda}$, such that the mass of $F$ is $\delta$. 
  If there exist at least two such geodesics,
  then this set is the set of parallel geodesics contained in a flat strip, of width $L>0$, hence it is homeomorphic to $[0,L]$. Hence, the set $F$ 
  is homeomorphic to a closed subset of $[0,L]$. Consequently, the preimage of $\nu$ by $\varphi_*$ is homeomorphic to the set of Borel measures on
$[0,L]$ of mass $\delta$.      

\medskip
  
 The spaces $\mathcal{M}\mathcal{L}_h(\Sigma)$ and $\mathcal{M}\mathcal{L}_p(\Sigma)$ are respectively homeomorphic to $\M_{\grperevet}([\Gmr])$ and to $\M_{\grperevet}([\Gqr])$, hence $\varphi_*$ 
  defines a continuous map $\psi:\mathcal{M}\mathcal{L}_p(\Sigma)\to\mathcal{M}\mathcal{L}_h(\Sigma)$ which is proper and surjective. The group $\RR^{+*}$ acts on these two spaces by multiplication
  of the measures. We denote by  
  $\mathcal{P}\mathcal{M}\mathcal{L}_p(\Sigma)$ and $\mathcal{P}\mathcal{M}\mathcal{L}_h(\Sigma)$
  the quotient spaces for these actions. Since $\psi$ is equivariant by these actions, it defines a continuous map $\overline{\psi}:
  \mathcal{P}\mathcal{M}\mathcal{L}_p(\Sigma)\to
  \mathcal{P}\mathcal{M}\mathcal{L}_h(\Sigma)$ which is surjective and proper. We deduce from this the following lemmas.
  
  \blemm The space $\mathcal{P}\mathcal{M}\mathcal{L}_p(\Sigma)$ is compact.
  \elemm
  
  \dem The space $\mathcal{P}\mathcal{M}\mathcal{L}_h(\Sigma)$ is compact (see \cite[Cor.~5~and~Prop.~17]{Bonahon88}) and $\overline{\psi}$ is proper.\cqfd
  
  \medskip
  
  If $\Sigma$ is compact, a {\it measured cylinder lamination} is a measured flat lamination having a unique component which is a cylinder component (see 
  \cite[§6]{Morzy1}).
  
  \blemm If $\Sigma$ is compact, the set of measured cylinder laminations having finitely many leaves is dense in $\mathcal{M}\mathcal{L}_p(\Sigma)$. In particular, 
  $\mathcal{M}\mathcal{L}_p(\Sigma)$ is separable.
  \elemm
  
  \dem The set of simple closed local geodesics endowed with a transverse measure which is a Dirac measure of positive mass is dense in $\mathcal{M}\mathcal{L}_h(\Sigma)$ (see
  \cite[Prop.~15]{Bonahon97}), and its preimage by $\varphi_*$ is the set of measured cylinder laminations. Since $\varphi_*$ is continuous, 
  this set is dense in $\mathcal{M}\mathcal{L}_p(\Sigma)$. If $(\Lambda,\mu)$ is a measured cylinder lamination whose support is not reduced to a single leaf, we denote by 
  $\alpha:[0,T]\to C$ ($T>0$) a geodesic arc such that the interior of the maximal flat cylinder $C$ that contains the support of $\Lambda$, endowed with the induced distance,
  is isometric to $\alpha(]0,T[)\times\SSS^1$.
  Then, the set of local geodesics contained in $C$ and parallel to the boundary of $C$ is homeomorphic to  
  $[0,T]$, and since the set of Radon measures with finite support on $[0,T]$ is dense in the space of Radon measures on $[0,T]$, there exists a sequence 
  $(\mu_{n,\alpha})_{n\in\NN}$ of Radon measures with finite support on $G(\alpha)$ such that, for all $n$, every leaf of the support 
  of $\mu_{n,\alpha}$ is parallel to the boundary of $C$ and $\mu_{n,\alpha}\overset{*}{\rightharpoonup}\mu_\alpha$.
  Moreover, $\Lambda_n=[\Supp(\mu_{n,\alpha})]$ is a flat lamination and $\mu_{n,\alpha}$ defines a transverse measure on $\Lambda_n$
  such that the sequence $(\Lambda_n,\mu_n)_{n\in\NN}$ converges to $(\Lambda,\mu)$.\cqfd
  
  \medskip
  
  The maps $\psi$ and $\overline{\psi}$ are not injective. Assume that $\Sigma$ is compact with genus  $g\in\NN$ and $b\in\NN$ boundary components.
  Considering the lack of injectivity of $\varphi_*$, the preimage of a measured hyperbolic lamination having no closed leaf consists in a unique measured flat lamination
  having 
  no periodic leaf. However, if $(\Lambda_m,\mu_m)$ is a closed leaf $\lambda$ endowed with a transverse measure which, for every arc $\alpha$ such that $G(\alpha)$ 
  contains $\lambda$, is a Dirac measure at $\lambda$ of mass $\delta>0$,
  the preimage of 
  $(\Lambda_m,\mu_m)$ by $\psi$ is the set of measured cylinder laminations whose supports are closed sets $F$ of leaves that are freely homotopic to $\lambda$. 
  If the set of local geodesics of $\srfce$ that are freely homotopic to $\lambda$ contains at least two elements, 
  it foliates a maximal flat cylinder. Then, this set is homeomorphic to $[0,L]$, with $L$ is the height of the flat cylinder, so $F$ is homeomorphic to a closed subset
  of $[0,L]$. Hence, the preimage of $(\Lambda_m,\mu_m)$ by $\psi$ is homeomorphic to the set of Borel measures on $[0,L]$ of mass $\delta$. Since $\varphi_*$ is 
equivariant for the addition of measures, we see that if a measured hyperbolic lamination has some closed leaves $\lambda_1,\dots,\lambda_p$, of respective masses 
 $\delta_1,\dots,\delta_p$  ($p$ is always at most $3g-3+b$), then the preimage of $(\Lambda_m,\mu_m)$ by $\psi$ is homeomorphic to the Cartesian product of the sets of
 Borel measures
 on $[0,L_i]$, $1\leqslant i\leqslant p$, where $L_i$ is the height of the maximal flat cylinder, union of the local geodesics of $\srfce$ freely homotopic to 
 $\lambda_i$, whose total mass is $\delta_i$.
  
  Since  $\Sigma$ is compact, the projectified space $\mathcal{P}\mathcal{M}\mathcal{L}_h(\Sigma)$ is homeomorphic to the sphere $\SSS^{6g-6+2b}$ 
  (see \cite[Th.~17]{Bonahon97}). 
  If the support of the measured hyperbolic laminations in the equivalence class of $x\in\mathcal{P}\mathcal{M}\mathcal{L}_h(\Sigma)$ has no 
  closed leaf, then the preimage of $\{x\}$
  by $\overline{\psi}$ is a single point. However, if the closed leaves $\lambda_1,\dots\lambda_p$ belong to the support of the measured hyperbolic 
  laminations in the equivalence
  class of $x$, and if $L_1,\dots,L_p$ are the heights of the maximal flat cylinders, unions of the images of the periodic local geodesics of $\srfce$ 
  that are freely homotopic 
  to $\lambda_1,\dots,\lambda_p$ ($L_i=0$ if there is only one periodic local geodesic that is freely homotopic to $\lambda_i$), then the preimage of 
  $\{x\}$ by
  $\overline{\psi}$ is homeomorphic to the Cartesian product with $p$ terms of the sets of Borel measures on  $[0,L_i]$ of total masses at most $1$.
  In particular, since the set of measured hyperbolic laminations whose support contains a closed leaf is projectively dense in the space of measured hyperbolic laminations,
  the subset of $\mathcal{P}\mathcal{M}\mathcal{L}_h(\Sigma)$ of points whose preimages by $\overline{\psi}$ are not a single point is dense in 
  $\mathcal{P}\mathcal{M}\mathcal{L}_h(\Sigma)$.

\subsection{Intersection number.}\label{nombreintersection}

If $g:\RR\to\widetilde{\Sigma}$ is a geodesic of $\revet$ or of $\revetm$, the orientation of $\widetilde{\Sigma}$ allows 
to define the two sides (arbitrarily) $+$ and $-$ of $g(\RR)$ (except
if $g(\RR)$ is a boundary component of $\partial\widetilde{\Sigma}$, in which case $g(\RR)$ has one side).
  We denote by $C_+(g)$ and $C_-(g)$ the two {\it complementary sides} which are the unions of $g(\RR)$ with the connected components of
  $\widetilde{\Sigma}-g(\RR)$ corresponding to one and the other sides of $g(\RR)$ 
  ($\widetilde{\Sigma}-g(\RR)$ may have more than two connected components if $g(\RR)$ is not disjoint from $\partial\widetilde{\Sigma}$).
  If $\widetilde{\Lambda}$ is a geodesic
  lamination of $\revet$ or $\revetm$, and if $g$ is a leaf of $\widetilde{\Lambda}$, then the image of another leaf of
  $\widetilde{\Lambda}$
is contained in $C_+(g)$ or $C_-(g)$, since the leaf is not interlaced with $g$. A leaf $g$ {\it separates} two other 
leaves if 
the image of one is contained in $C_+(g)$, and the image of the otherone in $C_-(g)$. Let $g_0$ and 
$g_1$ be two leaves of $\widetilde{\Lambda}$. We denote by $C_i$ the complementary side of $g_i(\RR)$ that contains $g_{i+1}(\RR)$ (with 
$i\in\ZZ/2\ZZ$). We define $C(g_0,g_1)=C_0\cap C_1$. Let $c$ be a geodesic segment joining their images (if the images are not disjoint, the segment $c$ may be a 
point). A leaf of
$\widetilde{\Lambda}$ {\it intersects $c$ non trivially } if it is contained in $C(g_0,g_1)$ and if it intersects both complementary components of the image of $c$ in
$C(g_0,g_1)$, and we denote by $B(g_0,g_1)=B_{\widetilde{\Lambda}}(g_0,g_1)$ the set of leaves of $\widetilde{\Lambda}$ intersecting
$c$ non trivially.

\blemm The compact set $B(g_0,g_1)$ does not depend on the choice of $c$.
\elemm

\dem The set of leaves whose images are contained in $C_0\cap C_1$ is compact and the requirement to intersect $c$ non trivially is closed in this set, hence  
$B(g_0,g_1)$ is compact.
Let $c'$ be another geodesic segment joining the images of $g_0$ and $g_1$. Since both $c$ and $c'$ separate $C_0\cap C_1$ in two connected
components and since the intersection of the images of two geodesic segments of $\revet$ (or between the image of a geodesic and a point) is connected,
every leaf that intersects $c$ non trivially also intersects  $c'$ non trivially, and conversely.\cqfd

  \medskip

The set of free homotopy classes of simple closed curves on $\Sigma$ endowed with transverse measures which are Dirac measures of positive masses, embeds into 
$\mathcal{M}\mathcal{L}_h(\Sigma)$, and the intersection number on this set can be extended, in a unique way,  to a continuous map 
$i:\mathcal{M}\mathcal{L}_h(\Sigma)\times\mathcal{M}\mathcal{L}_h(\Sigma)\to\RR^+$ (see \cite[Prop.~3]{Bonahon88}).
  According to Lemma \ref{surjectionpropre}, the map $\varphi_*$ defines a map $\psi:\mathcal{M}\mathcal{L}_p(\Sigma)\to\mathcal{M}\mathcal{L}_h(\Sigma)$.   
   Let $\alpha$ be a non trivial free homotopy class of closed curves, let $(\Lq,\mu_{[q]})$ be a measured flat lamination and let
  $\nu_{\mu_{[q]}}\in\M_{\grperevet}([\Gqr])$ be the measure defined by $(\Lq,\mu_{[q]})$ on $[G_{[\widetilde{q}]}]$ 
 (see Corollary \ref{correspondance}). We define the intersection number between $(\Lq,\mu_{[q]})$ and $\alpha$ by $$i_{[q]}(\mu_{[q]},\alpha)=
 i(\psi(\Lq,\mu_{[q]}),\alpha).$$ If $k\in\NN$, we have $i_{[q]}(\mu_{[q]},\alpha^k)=k\, i_{[q]}(\mu_{[q]},\alpha)$. Hence, we assume that $\alpha$ is primitive 
 (i.e. if there exists a free homotopy class $\alpha_0$ such that $\alpha=\alpha_0^k$, 
 then $k=\pm1$). We denote by $\alpha_{[q]}$ a periodic local flat geodesic  in the class of $\alpha$, and by $\widetilde{\alpha}_{[q]}$ a lift of
 $\alpha_{[q]}$ in $\widetilde{\Sigma}$. Let $\gamma\in\grperevet-\{e\}$ be one of the two primitive hyperbolic elements of $\grperevet$ whose translation axis is 
 $\widetilde{\alpha}_{[q]}(\RR)$, and let $(\widetilde{\Lambda}_{[q]},\mutilde_{[q]})$ be the preimage of $(\Lambda_{[q]},\mu_{[q]})$ in $\widetilde{\Sigma}$.
  
 \blemm\label{masseintersection} Let $\widetilde{\ell}$ be a leaf of $\Lqr$ which is interlaced with $\widetilde{\alpha}_{[q]}$. 
 The number $i_{[q]}(\mu_{[q]},\alpha)$ is equal to $\frac{1}{2}\nu_{\mu_{[q]}}(B_{\widetilde{\Lambda}}
 (\widetilde{\ell},\gamma\widetilde{\ell})-\gamma\widetilde{\ell})$. If there exists no such leaf, then $i_{[q]}(\mu_{[q]},\alpha)=0$.
 \elemm
 
 \dem The number $i(\psi(\Lq,\mu_{[q]}),\alpha)$ is equal to $\frac{1}{2}\varphi_*\nu_{\mu_{[q]}}(F_m)$, where $F_m$  is the set of leaves of 
 the measured hyperbolic 
 lamination $(\Lmr,\widetilde{\mu}_m)$ defined by $\varphi_*\nu_{\mu_{[q]}}$ which intersect a segment $I=[a,\gamma a[$ transversally, with
 $a\in\widetilde{\alpha}_m(\RR)$, where $\widetilde{\alpha}_m(\RR)$ is the translation axis
of $\gamma$ in $\revetm$, which is a fundamental domain of $\widetilde{\alpha}_m(\RR)$ for the action by translations of $\gamma^\ZZ$ (see \cite[Prop.~3]{Bonahon88}).
 Since the choice of $a$ is arbitrary, if $F_m$ is not empty, we may assume that $a$ is an intersection point between a leaf $\widetilde{\lambda}$ of $\Lmr$ and 
 $\widetilde{\alpha}_m(\RR)$. Hence, $F_m=B_{\Lmr}(\widetilde{\lambda},\gamma\widetilde{\lambda})-\gamma\widetilde{\lambda}$. 
 Moreover $\varphi_*\nu_{\mu_{[q]}}(F_m)=\nu_{\mu_{[q]}}(\varphi^{-1}(F_m))$, and by definition of $\varphi$, $\varphi^{-1}(F_m)=
 B_{\Lqr}(\widetilde{\ell},\gamma\widetilde{\ell})-\gamma\widetilde{\ell}$, 
 with $\widetilde{\ell}\in\Lqr$ belonging to $\varphi^{-1}(\lambda)$ (since $\nu_{\mu_{[q]}}$ is $\gamma$-invariant, if several leaves belong to 
 $\varphi^{-1}(\lambda)$, we can choose $\widetilde{\ell}$ arbitrarily in this set). 
  Finally, if $F_m$ is empty, no leaf of $\Lqr$ is interlaced with $\widetilde{\alpha}_{[q]}$ and 
 $i(\psi(\Lq,\mu_{[q]}),\alpha)=0$.\cqfd
 
 \medskip

 \rem We could define the intersection number between a free homotopy class of closed curves with $(\Lambda_{[q]},\mu_{[q]})$ by the infimum of the masses given
 by the measured
 flat lamination to the closed curves that are piecewise transverse to the lamination, similarly to the intersection number with a measured foliation, 
 but this infimum would not be necessarly attained since the periodic local geodesics are generally not piecewise transverse to the lamination.
 
 Furthermore, in the case of  compact surfaces endowed with a half-translation structure, whose boundary is empty, contrarily to measured hyperbolic lamination (see 
 \cite[Th.~2]{Ota90}), the intersection numbers with the free homotopy classes of closed curves of $\Sigma$ do not separate the measured flat laminations, but only 
 their images in  
 $\mathcal{M}\mathcal{L}_h(\Sigma)$. In particular, the topology defined after Definition \ref{defmesuretransverse} is not equivalent to the one induced by the product topology on 
 $\RR^{{\cal{H}}}$, with ${{\cal{H}}}$ the set of free homotopy classes of closed curves, on the image of $\mathcal{M}\mathcal{L}_p(\Sigma)$ by the map 
 $(\Lambda_{[q]},\mu_{[q]})\mapsto (i(\mu_{[q]},\alpha))_{\alpha\in{{\cal{H}}}}$. 
\section{Tree associated to a measured flat lamination.}\label{arbredual}

Let $\srfce$ be a compact, connected, orientable surface with (possibly empty) boundary, endowed with a half-translation structure, such that $\chi(\Sigma)<0$. Let 
$p:\revet\to\srfce$ be a locally isometric universal cover and let $(\Lambda,\mu)$ be a measured flat lamination on $\srfce$. We denote by 
$(\widetilde{\Lambda},\widetilde{\mu})$ its preimage in $\revet$ and by $\nu_{\widetilde{\mu}}$ the Borel measure that it defines on $\Gqr$ (see  Lemma \ref{nu}).
We first assume that $\nu_{\widetilde{\mu}}$ has no atom. 

\medskip

If $\{\widetilde{\ell}_1,\widetilde{\ell}_2\}$ is a pair of leaves of $\widetilde{\Lambda}$, we define  
$\widetilde{d}_{\widetilde{\Lambda}}(\widetilde{\ell}_1,\widetilde{\ell}_2)=\frac{1}{2}\nu_{\widetilde{\mu}}(B(\widetilde{\ell}_1,\widetilde{\ell}_2))$
(see Section \ref{nombreintersection} for the definition of $B(\widetilde{\ell}_1,\widetilde{\ell}_2))$). Then 
$\widetilde{d}_{\widetilde{\Lambda}}(\widetilde{\ell}_1,\widetilde{\ell}_2)\geqslant 0$ and
$\widetilde{d}_{\widetilde{\Lambda}}(\widetilde{\ell}_1,\widetilde{\ell}_2)=\widetilde{d}_{\widetilde{\Lambda}}(\widetilde{\ell}_2,\widetilde{\ell}_1)$. Moreover, if 
$\widetilde{\ell}_1$, $\widetilde{\ell}_2$ and $\widetilde{\ell}_3$ are three leaves of $\widetilde{\Lambda}$ and $c_1$, $c_2$ and $c_3$ are some geodesic segments
joining respectively the images of $\widetilde{\ell}_1$ and $\widetilde{\ell}_2$, $\widetilde{\ell}_2$ and $\widetilde{\ell}_3$, and $\widetilde{\ell}_1$ and 
$\widetilde{\ell}_3$ (denoted by $c_{1,2,3}$ on the figure, cases $2$, $3$ and $4$), then either none of the three leaves separates the other ones (case $1$, 
see Section \ref{nombreintersection} for the definition
of a leaf separating to other leaves)
or one separates the other ones
(case $2$, $3$ and $4$). In case $1$, every leaf of $B(\widetilde{\ell}_1,\widetilde{\ell}_3)$ intersects $c_1$ or $c_2$ non trivially, hence 
$B(\widetilde{\ell}_1,\widetilde{\ell}_3)\subseteq B(\widetilde{\ell}_1,\widetilde{\ell}_2)\cup B(\widetilde{\ell}_2,\widetilde{\ell}_3)$. In case $2$, we have 
$B(\widetilde{\ell}_1,\widetilde{\ell}_3)=B(\widetilde{\ell}_1,\widetilde{\ell}_2)\cup B(\widetilde{\ell}_2,\widetilde{\ell}_3)$. In case $3$, we have 
$B(\widetilde{\ell}_1,\widetilde{\ell}_3)\subseteq B(\widetilde{\ell}_1,\widetilde{\ell}_2)$ and in case $4$, we have 
$B(\widetilde{\ell}_1,\widetilde{\ell}_3)\subseteq B(\widetilde{\ell}_2,\widetilde{\ell}_3)$. In any case, we have
$\widetilde{d}_{\widetilde{\Lambda}}(\widetilde{\ell}_1,\widetilde{\ell}_3)\leqslant\widetilde{d}_{\widetilde{\Lambda}}(\widetilde{\ell}_1,\widetilde{\ell}_2)
+\widetilde{d}_{\widetilde{\Lambda}}(\widetilde{\ell}_2,\widetilde{\ell}_3)$. Hence $\widetilde{d}_{\widetilde{\Lambda}}$ is a pseudo-distance on $\widetilde{\Lambda}$.
\begin{center}
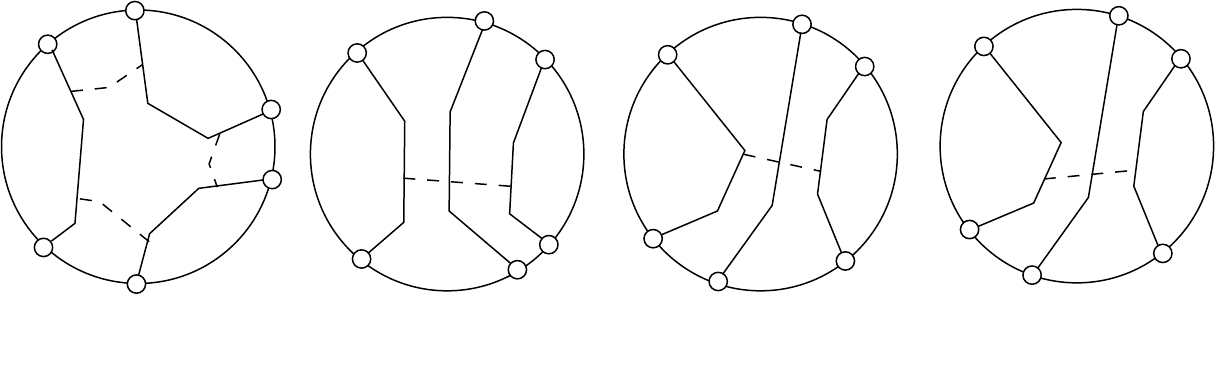
\end{center}

We denote by $(T,d_T)$ the  quotient metric space $(\widetilde{\Lambda},\widetilde{d}_{\widetilde{\Lambda}})/\sim$, where 
$\widetilde{\ell}\sim\widetilde{\ell}'$ if and only if $\widetilde{d}_{\widetilde{\Lambda}}(\widetilde{\ell},\widetilde{\ell}')=0$, and if
$F$ is a set of leaves of $\widetilde{\Lambda}$, we denote by $F^{\, T}$ its image by the quotient map. 

\brema\label{B(ell1ell2)} {\rm Let $\widetilde{\ell}_1$ and $\widetilde{\ell}_2$ be two distinct leaves of $\widetilde{\Lambda}$.
 Since $\nu_{\mutilde}$ has no atom and its support is $\widetilde{\Lambda}$, we have $\nu_{\mutilde}(B(\widetilde{\ell}_1,\widetilde{\ell}_2))=0$ if and only if 
$B(\widetilde{\ell}_1,\widetilde{\ell}_2)=\{\widetilde{\ell}_1,\widetilde{\ell}_2\}$, and the topology defined by the distance $d_T$ is equivalent to the quotient 
topology
of the topology induced by the geodesic topology on $\widetilde{\Lambda}$, by the equivalence relation $\widetilde{\ell}_1\;{\cal R}\;\widetilde{\ell}_2$ if and only if 
$B(\widetilde{\ell}_1,\widetilde{\ell}_2)=\{\widetilde{\ell}_1,\widetilde{\ell}_2\}$.

Finally, since the image of a leaf has two complementary sides (except if its image is a boundary component), and the image of each of the other leaves is contained in one of them,
the relation defined on $B(\widetilde{\ell}_1,\widetilde{\ell}_2)$ by $\widetilde{\ell}\preceq\widetilde{\ell}'$ if $\widetilde{\ell}$ belongs to
$B(\widetilde{\ell}_1,\widetilde{\ell}')$ is a total order which is compatible with ${\cal R}$, hence defines a total order on 
$B(\widetilde{\ell}_1,\widetilde{\ell}_2)^{\, T}$.}
\erema

\blemm\label{arbredualaunelamination} The metric space $(T,d_T)$ is an $\RR$-tree.
\elemm

\dem Let $\widetilde{\ell}_1$ and $\widetilde{\ell}_2$ be two leaves of $\widetilde{\Lambda}$. The map
$f:B(\widetilde{\ell}_1,\widetilde{\ell}_2)\to\RR^+$ 
defined by $f(\widetilde{\ell})=\widetilde{d}_{\widetilde{\Lambda}}(\widetilde{\ell}_1,\widetilde{\ell})$ is nondecreasing (with respect to $\preceq$) and continuous since 
$\nu_{\widetilde{\mu}}$ has no atom. 
Moreover, it is compatible with the equivalence relation ${\cal R}$ and defines an increasing continuous map $\overline{f}:B(\widetilde{\ell}_1,\widetilde{\ell}_2)^{\, T}
\to\RR^+$.
Since $B(\widetilde{\ell}_1,\widetilde{\ell}_2)^{\, T}$ is compact, it is a homeomorphism onto its image. Assume that its image is not an interval. Since it is a compact subset 
of $\RR$, if $U$ is a bounded  complementary component of ${f}(B(\widetilde{\ell}_1,\widetilde{\ell}_2))$ in $\RR^+$, then its closure is an interval $[a,b]$ with $a<b$.
Let $\widetilde{\ell}_a$ and $\widetilde{\ell}_b$ be some leaves of ${f}^{-1}(a)$ and ${f}^{-1}(b)$. If there exists a leaf $\widetilde{\ell}\in
B(\widetilde{\ell}_a,\widetilde{\ell}_b)-\{\widetilde{\ell}_a,\widetilde{\ell}_b\}$, since $\nu_{\mutilde}$ has no atom and its support is $\widetilde{\Lambda}$,
   we have $a<\overline{f}(\widetilde{\ell}^{\, T})<b$. Hence  
$B(\widetilde{\ell}_a,\widetilde{\ell}_b)=\{\widetilde{\ell}_a,\widetilde{\ell}_b\}$, which is impossible since then we would have   
$\widetilde{d}(\widetilde{\ell}_{a},\widetilde{\ell}_{b})=0$ and thus $a=b$. Hence, the image of $\overline{f}$ is the interval 
$[0,d_T(\widetilde{\ell}^{\, T}$,$\widetilde{\ell}'^{\, T})]$, and 
$\overline{f}^{-1}:[0,d_T(\widetilde{\ell}^{\, T}$, $\widetilde{\ell}'^{\, T})]\to T$ is a geodesic segment between $\widetilde{\ell}^{\, T}$ and 
$\widetilde{\ell}'^{\, T}$. Moreover, up to reparametrization, it is the unique arc joining $\widetilde{\ell}^{\, T}$ to $\widetilde{\ell}'^{\, T}$. Indeed, let $g:[0,1]\to T$ be another
arc joining $\widetilde{\ell}^{\, T}$ and $\widetilde{\ell}'^{\, T}$. If a leaf $\widetilde{\ell}_0$ belongs to $B(\widetilde{\ell},\widetilde{\ell}')$ then it separates 
$\widetilde{\ell}$ and $\widetilde{\ell}'$ (in the sense of Section \ref{nombreintersection}) and since $g$ is continuous for the quotient topology of the geodesic topology by
the equivalence relation ${\cal R}$, the point $\widetilde{\ell}_0^{\, T}$ belongs to the image of $g$. Hence $B(\widetilde{\ell},\widetilde{\ell}')^{\, T}$ is contained in the image of $g$.

Assume that there exists an element $x=g(t)$, with $t\in]0,1[$, in the image of $g$ that does not belong to
$B(\widetilde{\ell},\widetilde{\ell}')^{\, T}$ and let $\widetilde{\ell}_x$ be a leaf representing $x$. Assume for a contradiction that $\widetilde{\ell}_x(\RR)$ is
contained in $C(\widetilde{\ell},\widetilde{\ell}')$ (case $1$ below, see Section \ref{nombreintersection} for the definition of $C(\widetilde{\ell},\widetilde{\ell}')$).
Then $B(\widetilde{\ell},\widetilde{\ell}_x)$ is the union of the compact sets 
$B(\widetilde{\ell},\widetilde{\ell}_x)\cap B(\widetilde{\ell},\widetilde{\ell}')$ and $B(\widetilde{\ell},\widetilde{\ell}_x)\cap B(\widetilde{\ell}_x,\widetilde{\ell}')$, 
whose intersection is $\{\widetilde{\ell}_x\}$.
Let $\widetilde{\ell}_y$ be the closest element to $\widetilde{\ell}$ in the compact set $(B(\widetilde{\ell},\widetilde{\ell}_x)-B(\widetilde{\ell},\widetilde{\ell}')\cap 
B(\widetilde{\ell},\widetilde{\ell}_x))\cup\{\widetilde{\ell}_x\}$ (for $\preceq$). By assumption on
$\widetilde{\ell}_x$, we see that $\widetilde{\ell}_y$ is also the closest element to $\widetilde{\ell}'$ in  
$B(\widetilde{\ell}_x,\widetilde{\ell}')-B(\widetilde{\ell},\widetilde{\ell}')\cap B(\widetilde{\ell}',\widetilde{\ell}_x)$.
If $\widetilde{\ell}_y\not=\widetilde{\ell}_x$, the leaf $\widetilde{\ell}_y$ separates $\widetilde{\ell}$ from $\widetilde{\ell}_x$ and $\widetilde{\ell}_x$ from 
$\widetilde{\ell}'$. Since $g$ is continuous,
the element $\widetilde{\ell}_y^{\, T}$ belongs to $g([0,t])$ and to $g([t,1])$.  
If $B(\widetilde{\ell}_y,\widetilde{\ell}_x)\not=\{\widetilde{\ell}_y,\widetilde{\ell}_x\}$, then $\widetilde{\ell}_x^{\, T}\not=\widetilde{\ell}_y^{\, T}$,
and $\widetilde{\ell}_y^{\, T}$ would belong
to $g([0,t[)$ and to $g(]t,1])$, thus $g$ would not be injective. If $B(\widetilde{\ell}_y,\widetilde{\ell}_x)=\{\widetilde{\ell}_y,\widetilde{\ell}_x\}$, we denote by 
$\widetilde{\ell}_z$ the closest element to $\widetilde{\ell}_y$ in the compact set 
$ B(\widetilde{\ell},\widetilde{\ell}_y)\cap B(\widetilde{\ell},\widetilde{\ell}')$. By definition of $\widetilde{\ell}_y$,
we have $B(\widetilde{\ell}_z,\widetilde{\ell}_y)=\{\widetilde{\ell}_z,\widetilde{\ell}_y\}$, thus $\widetilde{\ell}_y^{\, T}=\widetilde{\ell}_z^{\, T}$, and since
$\widetilde{\ell}_x^{\, T}=\widetilde{\ell}_y^{\, T}$, we have $\widetilde{\ell}_x^{\, T}\in B(\widetilde{\ell},\widetilde{\ell}')^{\, T}$, which is a contradiction.

\medskip

Hence $\widetilde{\ell}_x(\RR)$ is contained in the complementary side of $\widetilde{\ell}(\RR)$ that does not contain
$\widetilde{\ell}'(\RR)$, or the opposite (cases $2$ and $3$ below). However, since $g$ is continuous and since $\widetilde{\ell}$ separates $\widetilde{\ell}_x$ and
$\widetilde{\ell}'$ (or $\widetilde{\ell}'$ separates $\widetilde{\ell}_x$ and $\widetilde{\ell}$), there would then exist $t\in\mathopen{]}0,1[$ such that 
$\widetilde{\ell}$ (or $\widetilde{\ell}'$) represents $g(t)$,
et $g$ would not be injective.
\begin{center}
 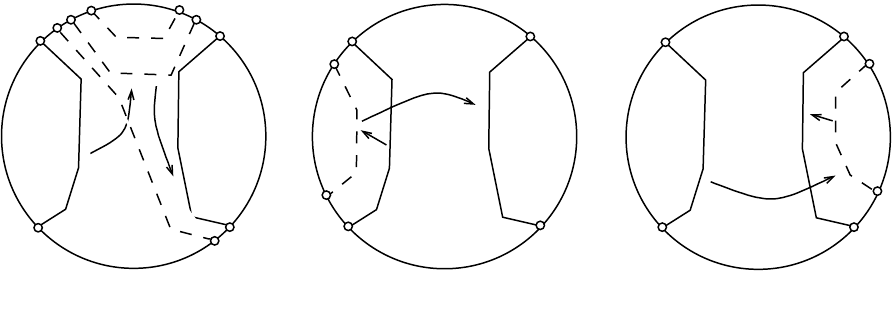
\end{center}

Hence, if $g$ is an arc between $\widetilde{\ell}^{\, T}$ and $\widetilde{\ell}'^{\, T}$, it has the same image as $\overline{f}^{-1}$. Thus, up to reparametrization, the 
unique arc between $\widetilde{\ell}^{\, T}$ and $\widetilde{\ell}'^{\, T}$ is $\overline{f}^{-1}$ which is isometric to 
$[0,d_T(\widetilde{\ell}^{\, T}$,$\widetilde{\ell}'^{\, T})]$ by construction. Since it is true for all pair of leaves, the metric space 
$(T,d_T)$ is an $\RR$-tree.\cqfd

\medskip

Assume that the measure $\nu_{\widetilde{\mu}}$ defined by $(\widetilde{\Lambda},\widetilde{\mu})$ has an atom $\widetilde{\ell}$.
We replace $\widetilde{\ell}(\RR)$ by a flat strip of width $\nu_{\widetilde{\mu}}(\widetilde{\ell})$ by gluing isometrically each of the complementary sides of 
$\widetilde{\ell}(\RR)$ on the boundary components of 
$\operatorname{FS}(\widetilde{\ell})=\RR\times[0,\nu_{\widetilde{\mu}}(\widetilde{\ell})]$, endowed with the Euclidean distance.
By doing the same  thing for every atom of $\nu_{\widetilde{\mu}}$, we get
a surface $\widetilde{\Sigma}'$ endowed with a complete $\CAT(0)$ metric $\widetilde{d}'$ and the isometric action of the covering group $\Gamma_{\widetilde{\Sigma}}$
on $\widetilde{\Sigma}$ minus the images of the atoms of
$\nu_{\widetilde{\mu}}$ extends in a unique way to an isometric action on $(\widetilde{\Sigma}',\widetilde{d}')$. Note that $(\widetilde{\Sigma},
\widetilde{d}')$ is a complete $\CAT(0)$ metric space whose boundary at infinity is endowed with a (total) cyclic order, compatible
with the covering group action, hence it enters in the generalized framework of the geodesic
laminations introduced in \cite{Morzy1}. However, the distance $\widetilde{d}'$ does not necessarily come from a half-translation structure, since the angles of the
conical singular points belonging to the boundary of an added flat strip may not be multiples of $\pi$.

Let $\widetilde{\ell}$ be an atom of $\nu_{\widetilde{\mu}}$ and let $F_{\widetilde{\ell}}$ be the maximal set of geodesics whose images are contained
in the corresponding flat strip $\BP(\widetilde{\ell})$ of $(\widetilde{\Sigma}',\widetilde{d}')$, that are parallel to its boundary. Let $\alpha$ be a geodesic
segment of $\BP(\widetilde{\ell})$ that orthogonally joins its boundary components. Then, the map
$r:F_{\widetilde{\ell}}\to\Image(\alpha)$ defined by 
$r(g)=g(\RR)\cap\Image(\alpha)$ is a homeomorphism. Hence, we can endow $F_{\widetilde{\ell}}$ with the measure $\nu_{\widetilde{\ell}}=
(r^{-1})_*dx_ \alpha$, where 	
$dx_\alpha$ is the Lebesgue measure on $\Image(\alpha)$, of total mass $\nu_{\widetilde{\mu}}(\widetilde{\ell})$. 
The canonical isometric embedding of every connected component of $\widetilde{\Sigma}-
\bigcup_{\nu_{\mutilde}(\widetilde{\ell})>0}\widetilde{\ell}(\RR)$ into $\widetilde{\Sigma}'$ induces an embedding of $\widetilde{\Lambda}-
\bigcup_{\nu_{\mutilde}(\widetilde{\ell})>0}\widetilde{\ell}$ into the space $[\G_{\widetilde{d}'}]$ of geodesics of $(\widetilde{\Sigma}',\widetilde{d}')$ defined up to
changing origin. We denote by $\nu'_c$ the push forward of $\nu_{\mutilde|\widetilde{\Lambda}-
\bigcup_{\nu_{\mutilde}(\widetilde{\ell})>0}\widetilde{\ell}}$ by this embedding. Then, we define $\nu_{\mutilde}'=\nu'_c+\sum_{\nu_{\mutilde}
(\widetilde{\ell})>0}\nu_{\widetilde{\ell}}$.
The measure $\nu_{\mutilde}'$ is a $\grperevet$-invariant Radon measure on the set of geodesics (defined up to changing origin) of $(\widetilde{\Sigma}',\widetilde{d}')$. 
Its support
is the union of the image of $\widetilde{\Lambda}-\bigcup_{\nu_{\mutilde}(\widetilde{\ell})>0}\widetilde{\ell}$ by the canonical embedding, and of the sets 
$F_{\widetilde{\ell}}$, with $\nu(\widetilde{\ell})>0$, defined above, which is a geodesic lamination  $\widetilde{\Lambda}'$ of $(\widetilde{\Sigma}',\widetilde{d}')$. 

Since, by construction, the Radon measure $\nu'_{\mutilde}$ has no atom, we can define the tree associated to the measure $\nu'_{\mutilde}$, exactly as above, even if 
its support is not a flat lamination as defined in this paper, since it only matters than $(\widetilde{\Sigma}',\widetilde{d}')$ is complete, $\CAT(0)$, and the support $\widetilde{\Lambda}'$ of
$\nu'_{\mutilde}$ is a geodesic lamination in the generalized sense of \cite{Morzy1}. We will call {\it tree associated to $(\Lambda,\mu)$ } the tree associated to 
$\nu_{\mutilde}'$ defined in this way.

%
\section{Covering group action on the tree associated to a measured flat lamination.}\label{actiondegroupe}

In this section, we use the definitions and notations of Section \ref{arbredual} and we consider the canonical action of the covering group  
$\Gamma_{\widetilde{\Sigma}}$ on the tree $(T,d_T)$ associated to the preimage $(\widetilde{\Lambda},\widetilde{\mu})$ of the measured flat lamination  
$(\Lambda,\mu)$ in $\widetilde{\Sigma}$. We may assume that $\nu_{\widetilde{\mu}}$ has no atom, up to proceeding 
as in the last paragraph of Section \ref{arbredual}.
Since $\widetilde{\Lambda}$ is fixed, we will use the notation  
$B(\widetilde{\ell},\widetilde{\ell}')$ instead of $B_{\widetilde{\Lambda}}(\widetilde{\ell},\widetilde{\ell}')$ for every pair of leaves $\widetilde{\ell},\widetilde{\ell}'$ of
$\widetilde{\Lambda}$. 

\medskip

The covering group $\Gamma_{\widetilde{\Sigma}}$ acts on $\widetilde{\Sigma}$ by isometries. Hence, it defines an action on the set  
$[\G_{[\widetilde{q}]}]$ of geodesics of $\revet$  that are defined up to changing origin. Since $\widetilde{\Lambda}$ is $\grperevet$-invariant, this action defines an action on
$\widetilde{\Lambda}$. Since for every $\gamma\in\grperevet$ we have $\gamma_*\nu_{\widetilde{\mu}}=\nu_{\widetilde{\mu}}$ and 
$\gamma B(\widetilde{\ell},\widetilde{\ell}')=B(\gamma\widetilde{\ell},\gamma\widetilde{\ell}')$, for every pair of leaves $\widetilde{\ell},\widetilde{\ell}'$ of
$\widetilde{\Lambda}$, it defines an isometric action of $\grperevet$ on the  tree $(T,d_T)$ associated to $(\widetilde{\Lambda},\widetilde{\mu})$ defined in Lemma 
\ref{arbredualaunelamination}.

\blemm For every primitive element $\gamma\in\grperevet-\{e\}$, if $\widetilde{\alpha}_{\gamma}(\RR)$ is a tranlation axis of $\gamma$ in $\revet$ and $\alpha_\gamma$
is the projection of $\widetilde{\alpha}_\gamma$ in $\Sigma$, then the translation distance $\ell_T(\gamma)$ of $\gamma$ in $(T,d_T)$ 
is equal to $i_{[q]}(\mu,\alpha_\gamma)$. Moreover, if $\ell_T(\gamma)>0$, the translation axis of $\gamma$ is the image in $T$ of the set of leaves of
$\widetilde{\Lambda}$ which are interlaced with $\widetilde{\alpha}_\gamma$.  
\elemm

\dem {\bf Case $(1)$.~} Assume that $i_{[q]}(\mu,\alpha_\gamma)>0$. Then $\widetilde{\alpha}_\gamma$ is interlaced with at least one leaf $\widetilde{\ell}$ of  
$\widetilde{\Lambda}$. Hence, acccording to Lemma \ref{masseintersection} and since  
$\nu_{\widetilde{\mu}}(\{\gamma\widetilde{\ell}\})=0$, we have $i_{[q]}(\mu,\alpha_\gamma)=\frac{1}{2}\nu_{\mutilde}(B(\widetilde{\ell},\gamma\widetilde{\ell}))
=\widetilde{d}(\widetilde{\ell},\gamma\widetilde{\ell})$.

Moreover, the set $F$ of leaves of $\widetilde{\Lambda}$ whose images are contained in $C(\widetilde{\ell},\gamma\widetilde{\ell})$, minus $\gamma\widetilde{\ell}$,
is a fundamental domain of 
$\widetilde{\Lambda}$ for the action of $\gamma^\ZZ$ (see Section \ref{nombreintersection} for the definition of $C(\widetilde{\ell},\gamma\widetilde{\ell})$). 
  If $\widetilde{\ell}'$ is a leaf of $\widetilde{\Lambda}$, there exists a unique $n\in\ZZ$ such that 
$\gamma^n\widetilde{\ell}'$ belongs to $F$. 

\begin{center}
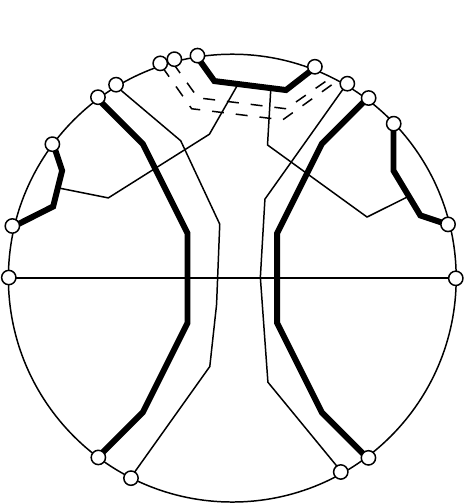
\end{center}

Let $\widetilde{\ell}'_0=\gamma^n\widetilde{\ell}'$, $\widetilde{\ell}_1'=\gamma\widetilde{\ell}'_0$ and 
$\widetilde{\ell}'_{-1}=\gamma^{-1}\widetilde{\ell}'_0$. Then $\widetilde{\ell}_{-1}'$ and $\widetilde{\ell}_1'$ do not belong to $F$. Since the leaves of 
$\widetilde{\Lambda}$ are pairwise non interlaced, the leaves of $B(\widetilde{\ell},\gamma\widetilde{\ell})$
belong to $B(\widetilde{\ell}'_{-1},\widetilde{\ell}'_{0})\cup B(\widetilde{\ell}'_{0},\widetilde{\ell}'_{1})$, and
$\gamma (B(\widetilde{\ell}'_{-1},\widetilde{\ell}'_0)\cap B(\widetilde{\ell},\gamma\widetilde{\ell})-\widetilde{\ell})\subseteq B(\widetilde{\ell}'_{0},
\widetilde{\ell}'_{1})-B(\widetilde{\ell}'_{0},\widetilde{\ell}'_{1})\cap B(\widetilde{\ell},\gamma\widetilde{\ell})$, 
since $F$ is a fundamental domain of $\widetilde{\Lambda}$ for the action of $\gamma^\ZZ$. Since $\nu_{\mutilde}$ is $\grperevet$-invariant, we have
\begin{align*}
\nu_{\mutilde}(B(\widetilde{\ell}',\gamma\widetilde{\ell}'))&=\nu_{\mutilde}(B(\widetilde{\ell}'_{0},\widetilde{\ell}'_1))\\
&\geqslant \nu_{\mutilde}(B(\widetilde{\ell},\gamma\widetilde{\ell})\cap B(\widetilde{\ell}'_0,\widetilde{\ell}'_1))+
\nu_{\mutilde}(\gamma(B(\widetilde{\ell},\gamma\widetilde{\ell})\cap B(\widetilde{\ell}'_{-1},\widetilde{\ell}'_0)))\\
&=\nu_{\mutilde}(B(\widetilde{\ell},\gamma\widetilde{\ell})\cap B(\widetilde{\ell}'_0,\widetilde{\ell}'_1))+
\nu_{\mutilde}(B(\widetilde{\ell},\gamma\widetilde{\ell})\cap B(\widetilde{\ell}'_{-1},\widetilde{\ell}'_0))\\
&\geqslant\nu_{\mutilde}(B(\widetilde{\ell},\gamma\widetilde{\ell}))\mbox{ since } B(\widetilde{\ell},\gamma\widetilde{\ell})\subseteq 
B(\widetilde{\ell}'_{-1},\widetilde{\ell}'_0)\cup B(\widetilde{\ell}'_0,\widetilde{\ell}'_1)\\
&=\widetilde{d}(\widetilde{\ell},\gamma\widetilde{\ell})
\end{align*}

Hence $\ell_T(\gamma)=\widetilde{d}(\widetilde{\ell},\gamma\widetilde{\ell})=i_{[q]}(\mu,\alpha_\gamma)$. In particular, we have $\ell_T(\gamma)>0$  and 
the isometric action of $\gamma$ on $(T,d_T)$ is hyperbolic, and hence admits a tranlation axis.   
Moreover, if $\widetilde{\ell}'$ is not interlaced with  
$\widetilde{\alpha}_\gamma$, then neither are $\widetilde{\ell}'_{-1}$, $\widetilde{\ell}'_0$ and $\widetilde{\ell}_1$. However, we have seen that
$B(\widetilde{\ell},\gamma\widetilde{\ell})$ is contained in $B(\widetilde{\ell}'_{-1},\widetilde{\ell}'_0)\cup B(\widetilde{\ell}'_0,\widetilde{\ell}'_1)$ and that
$ B(\widetilde{\ell}'_0,\widetilde{\ell}'_1)$ contains the union of the sets $B(\widetilde{\ell},\gamma\widetilde{\ell})\cap 
B(\widetilde{\ell}'_0,\widetilde{\ell}'_1)$ and $\gamma(B(\widetilde{\ell},\gamma\widetilde{\ell})\cap 
B(\widetilde{\ell}'_{-1},\widetilde{\ell}'_0))$, and their intersection is $\{\gamma\widetilde{\ell}\}$.     
Assume that $\widetilde{d}(\widetilde{\ell}'_0,\widetilde{\ell}'_{1})$ is equal to $\widetilde{d}(\widetilde{\ell},\gamma\widetilde{\ell})$, i.e. that
$\nu_{\mutilde}(B(\widetilde{\ell}'_0,\widetilde{\ell}'_1))=\nu_{\mutilde}(B(\widetilde{\ell},\gamma\widetilde{\ell})\cap 
B(\widetilde{\ell}'_0,\widetilde{\ell}'_1))+\nu_{\mutilde}(\gamma(B(\widetilde{\ell},\gamma\widetilde{\ell})\cap 
B(\widetilde{\ell}'_{-1},\widetilde{\ell}'_0)))$. Since the support of $\nu_{\mutilde}$ is $\widetilde{\Lambda}$, there are two possibilities.

In the first possibility, the intersection 
$(F\cup\gamma\widetilde{\ell})\cap(B(\widetilde{\ell}'_{0},\widetilde{\ell}'_1)
\cup \gamma^{-1}B(\widetilde{\ell}'_0,\widetilde{\ell}'_1))
=(F\cup\gamma\widetilde{\ell})\cap(B(\widetilde{\ell}'_{-1},\widetilde{\ell}'_0)\cup B(\widetilde{\ell}'_0,\widetilde{\ell}'_1))$ is equal to 
$B(\widetilde{\ell},\gamma\widetilde{\ell})$. Thus, no leaf of $F$ is contained in the complementary side $A$ of $\widetilde{\ell}'_0(\RR)$ that contains  
$\widetilde{\alpha}_\gamma(\RR)$, if it is non interlaced with $\widetilde{\alpha}_\gamma$ (there is no leaf as represented by a dotted line  points on the above picture).
Consequently, there exists $\widetilde{\ell}"\in B(\widetilde{\ell},\gamma\widetilde{\ell})$ such that the set $B(\widetilde{\ell}'_0,\widetilde{\ell}")$ is reduced to
$\{\widetilde{\ell}'_0,\widetilde{\ell}"\}$ and $\widetilde{d}(\widetilde{\ell}",\widetilde{\ell}'_0)=0$. Thus $\widetilde{\ell}_0'^{\, T}=\widetilde{\ell}"`^{\, T}\in 
B(\widetilde{\ell},\gamma{\widetilde{\ell}})^T$.

In the second possibility, there is a minimum, denoted by $m$ in $B(\widetilde{\ell}'_0,\widetilde{\ell})$, for $\preceq$ defined in Remark \ref{B(ell1ell2)},
different from  $\widetilde{\ell}'_0$.
Then $\widetilde{d}(m,\widetilde{\ell}'_0)=0$ and since $m$ belongs to $B(\widetilde{\ell}'_0,\widetilde{\ell})\subset 
F\cap(B(\widetilde{\ell}'_{-1},\widetilde{\ell}_0)$, it is interlaced with $\widetilde{\alpha}_\gamma$ by assumption. Hence, the image of $\widetilde{\ell}'_0$ in $T$ 
belongs to the translation axis of $\gamma$ in $T$ and, by $\gamma$-invariance, so does the image of $\widetilde{\ell}$'. Hence, if $i_{[q]}(\mu,\alpha_\gamma)>0$, 
the translation distance $\ell_T(\gamma)$ of $\gamma$ is equal to $i_{[q]}(\mu,\alpha_\gamma)$ and the translation axis of $\gamma$ is the image in $T$
of the set of leaves of 
$\widetilde{\Lambda}$ which are interlaced with a translation axis of $\gamma$ in $\revet$.

{\bf Case $(2)$.~} Assume that $i_{[q]}(\mu,\alpha_\gamma)=0$, or equivalently that $\widetilde{\alpha}_\gamma$ is interlaced with no leaf of $\widetilde{\Lambda}$. 
If $\widetilde{\alpha}_\gamma$ has the same ordered pair of points at infinity as a leaf $\widetilde{\ell}$ of $\widetilde{\Lambda}$,
then $\gamma\widetilde{\ell}=\widetilde{\ell}$ and $\widetilde{\ell}^{\,T}$ is a fixed point of $\gamma$ in $T$.
Otherwise, we  define 
$(S,N)=(\widetilde{\alpha}_\gamma(-\infty),\widetilde{\alpha}_\gamma(+\infty))$. According to Lemma \ref{quedescomposnatesminimales},
no leaf of $\Lambda$ is positively periodic unless it is periodic, hence according to \cite[Lem.~4.13~et~4.14]{Morzy1}, neither $N$ nor $S$
is a point at infinity of any leaf of $\widetilde{\Lambda}$. We recall that the (total) cyclic order $o$ on $\partial_\infty\widetilde{\Sigma}$, 
defined by the orientation of $\widetilde{\Sigma}$, 
defines a total order $\leqslant$ on $\partial_\infty\widetilde{\Sigma}$ defined by $S<\eta$ for every  
$\eta\in\partial_\infty\widetilde{\Sigma}-\{S\}$ and $\eta_1\leqslant\eta_2$ if and only if 
$o(\eta_1,\eta_2,S)\in\{0,1\}$ for every $\eta_1,\eta_2\in\partial_\infty\widetilde{\Sigma}-\{S\}$ 
(see \cite[Rem.~2.9]{Morzy1} for the definition of $o$ and \cite[Déf.~2.23]{Wolf11} for the definition of $\leqslant$).

Let $\widetilde{\ell}$ be a leaf of $\widetilde{\Lambda}$ and $(a,b)=(\widetilde{\ell}(-\infty),\widetilde{\ell}(+\infty))$. 
Since no leaf of $\widetilde{\Lambda}$ is interlaced with $\widetilde{\alpha}_\gamma$, the union $\widetilde{\Lambda}\cup
\{\widetilde{\alpha}_\gamma,\widetilde{\alpha}_\gamma^{-1}\}$ is a flat lamination, and we can define $B(\widetilde{\ell},\widetilde{\alpha}_\gamma)$.
We replace $\widetilde{\ell}$ by the maximum of (the compact set) 
$B(\widetilde{\ell},\widetilde{\alpha}_\gamma)-\widetilde{\alpha}_\gamma$, for $\preceq$. Then $B(\widetilde{\ell},\widetilde{\alpha}_\gamma)=
\{\widetilde{\ell},\widetilde{\alpha}_\gamma\}$.
However, the action of $\gamma^\ZZ$ on $\partial_\infty\widetilde{\Sigma}$ has a South-North dynamic, whose fixed points are $N$ and $S$, with $N$ attractive and $S$ 
repulsive. Up to taking the opposite of the cyclic order on $\partial_\infty\widetilde{\Sigma}$, we can assume that the action of $\gamma^\ZZ$ on the $\gamma^\ZZ$-orbits of $a$
and $b$ is increasing. 

\medskip
\noindent
\begin{minipage}{10 cm}
Hence $a<\gamma a$, $b<\gamma b$ and since $\widetilde{\ell}$ and $\gamma\widetilde{\ell}$ are not interlaced, we also have $b\leqslant \gamma a\leqslant \gamma b$.
Assume that there exists a leaf $\widetilde{\ell}'\in B(\widetilde{\ell},\gamma\widetilde{\ell})
-\{\widetilde{\ell},\gamma\widetilde{\ell}\}$, whose ordered pair of points at infinity is $(a',b')$. Since $\widetilde{\ell}'$ is not interlaced with 
$\widetilde{\alpha}_\gamma$ and by assumption on $\widetilde{\ell}$, up to replacing $\widetilde{\ell}'$  by its opposite, we have 
$b\leqslant a'\leqslant \gamma a$ and 
$\gamma b\leqslant b'< N$, and $(a',b')\not=(\gamma a,\gamma b)$. But then $S<\gamma^{-1}a'\leqslant a$ and $b\leqslant\gamma^{-1}b'< N$, with
$(\gamma^{-1}a',\gamma^{-1}b')\not=(a,b)$, thus $\gamma^{-1}\widetilde{\ell}'\in B(\widetilde{\ell},\widetilde{\alpha}_\gamma)-\{\widetilde{\ell},
\widetilde{\alpha}_\gamma\}$, which is a contradiction to the assumption on $\widetilde{\ell}$. Hence 
$B(\widetilde{\ell},\gamma\widetilde{\ell})=\{\widetilde{\ell},\gamma\widetilde{\ell}\}$ and $\widetilde{d}(\widetilde{\ell},\gamma\widetilde{\ell})=0$.
Hence $\widetilde{\ell}^{\,T}$ is a fixed point of $\gamma$ in $T$.\cqfd
\end{minipage}
\begin{minipage}{4.9 cm}
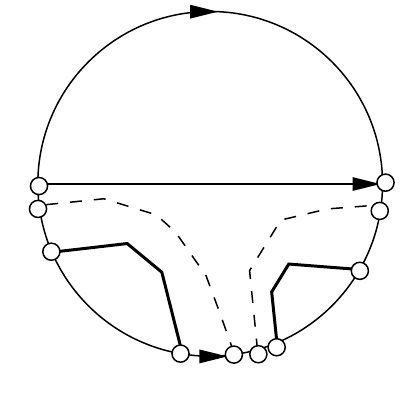
\end{minipage}

\section{Isometry between the tree associated to a measured flat lamination and the dual tree to the corresponding measured hyperbolic lamination.}\label{linkstrees}

 In this section \ref{linkstrees}, we use the same notation as in Section \ref{arbredual}. We denote by $m$ a hyperbolic metric with totally geodesic boundary on
 $\Sigma$ and by $\widetilde{m}$ the lifted hyperbolic metric on 
 $\widetilde{\Sigma}$.  
 We begin by recalling the definition of the dual tree to a measured hyperbolic lamination (see for example \cite[§1]{MorSha91}), with a new presentation that allows 
 to construct a $\grperevet$-equivariant isometry between the tree associated to a measured flat lamination and the dual tree to the corresponding measured hyperbolic 
 lamination. Let $(\Lambda_m,\mu_m)$ be a measured hyperbolic lamination of $\srfcem$ and let $(\Lmr,\mutilde_{{m}})$ be its preimage in $\widetilde{\Sigma}$. Then
 $(\Lmr,\mutilde_{{m}})$ is $\grperevet$-invariant and defines a measure 
$\nu_{\mutilde_m}\in\M_{\grperevet}([\Gmr])$ (see \cite[Prop.~17~p.~154]{Bonahon88}). If $\widetilde{\lambda}$ is an atom of $\nu_{\mutilde_m}$, we replace  
$\widetilde{\lambda}$ by a flat strip $\BP(\widetilde{\lambda})$ of width $\nu_{\mutilde_m}(\widetilde{\lambda})$ foliated by geodesic lines parallel to its boundary. 
Proceeding similarly for every
atom of $\nu_{\mutilde_m}$, we get a geodesic lamination $\widetilde{\Lambda}'$ on a surface endowed with a distance $(\widetilde{\Sigma}',
d')$, which is $\CAT(0)$, and the isometric action of $\grperevet$ on $\revetm$ extends uniquely
to an isometric action on $(\widetilde{\Sigma}',d')$. Then, it is the same framework as at the end of Section \ref{arbredual}. We can similarly defines the measure 
$\nu_{\widetilde{\lambda}}$ of total mass $\nu_{\widetilde{\mu}_m}(\widetilde{\lambda})$ on the parallel leaves foliating the flat strip $\BP(\widetilde{\lambda})$ 
associated to an atom
$\widetilde{\lambda}$ of $\nu_{\mutilde_m}$. The measure  
$\nu'$ which is equal to $\nu_{\mutilde_m}$ outside of the atoms and equal to $\nu_{\widetilde{\lambda}}$ on the set of leaves foliating the flat stip associated to 
$\widetilde{\lambda}$, for every atom $\widetilde{\lambda}$ of $\nu_{\mutilde_m}$, is a Radon measure on $[\G_{d'}]$ which is $\grperevet$-invariant, atomless,
 and whose support is equal to $\widetilde{\Lambda}'$.

We also define a pseudo-distance
$\widetilde{d}_{\widetilde{\Lambda}'}$ on $\widetilde{\Lambda}'$ by
$\widetilde{d}_{\widetilde{\Lambda}'}(\widetilde{\lambda}_0,\widetilde{\lambda}_1)=\frac{1}{2}\nu'(B_{\widetilde{\Lambda}'}(\widetilde{\lambda}_0,\widetilde{\lambda}_1))$
for every pair of leaves $\widetilde{\lambda}_0$, $\widetilde{\lambda}_1$ of $\widetilde{\Lambda}'$ (see Section \ref{nombreintersection} for the definition of 
$B_{\widetilde{\Lambda}'}(\widetilde{\lambda}_0,\widetilde{\lambda}_1)$), and the quotient  
of $(\widetilde{\Lambda}',d_{\widetilde{\Lambda}'})$ by the equivalence relation $\widetilde{\lambda}_0\sim\widetilde{\lambda}_1$ if and only if 
$d_{\widetilde{\Lambda}'}(\widetilde{\lambda}_0,\widetilde{\lambda}_1)=0$ (or equivalently $B(\widetilde{\lambda}_0,\widetilde{\lambda}_1)=
\{\widetilde{\lambda}_0,\widetilde{\lambda}_1\}$) is an $\RR$-tree $(T,d_T)$ called the {\it tree dual} to $(\Lmr,\mutilde_m)$. 
For every 
$\gamma\in\grperevet$, we have
$\gamma_*\nu'=\nu'$ and 
$\gamma B(\widetilde{\lambda}_0,\widetilde{\lambda}_1)=B(\gamma\widetilde{\lambda}_0,\gamma\widetilde{\lambda}_1)$ for every pair of leaves $\widetilde{\lambda}_0$ and
$\widetilde{\lambda}_1$ of $\widetilde{\Lambda}'$. Hence, the action of $\grperevet$ on $\Lmr$ defines an isometric action on $(T,d_T)$. It is easy to check (compare for 
instance with \cite{LevPau97}) that there exists 
a $\grperevet$-equivariant isometry from  the dual tree constructed in this way, and the one constructed for example in \cite[§1]{MorSha91}, by identification of the 
leaves of $\widetilde{\Lambda}'$ with the complementary connected component of the support of $\widetilde{\Lambda}'$ that they bound.

\medskip

Let $(\Lqr,\mutilde_{[\widetilde{q}]})$ be a measured flat lamination on $\revet$, let $\nu_{\mutilde_{[{q}]}}$ be its associated measure on $[\G_{[\widetilde{q}]}]$
and let $\nu_{\mutilde_m}$ be its image by $\varphi_*$ (see Lemma \ref{surjectionpropre}). We denote by $(\Lmr,\mutilde_{m})$ the measured hyperbolic lamination 
defined by $\nu_{\mutilde_m}$, and by $\widetilde{\Lambda}'$ and $\nu'$ the geodesic lamination on $(\widetilde{\Sigma}',d')$ and the Radon measure on $[\G_{d'}]$
defined by $(\Lmr,\mutilde_m)$ as above. We assume (up to proceeding as in the last paragraph of Section \ref{arbredual}) that $\nu_{\mutilde_{[{q}]}}$ has no atom.

If $\widetilde{\lambda}$ is an atom of $\nu_{\mutilde_m}$, and if $F_{\widetilde{\lambda}}$ is the set of leaves of $\Lqr$ to which corresponds
$\widetilde{\lambda}$ 
(see \cite[§4.2]{Morzy1}), there exist maximal flat strips in $\revet$ and in $(\widetilde{\Sigma}',d')$ that contain respectively $F_{\widetilde{\lambda}}$ and the set 
$F'_{\widetilde{\lambda}}$ of leaves of $\widetilde{\Lambda}'$ corresponding to $\widetilde{\lambda}$ in the above construction. We denote by $\widetilde{\ell}_0$ and 
$\widetilde{\ell}_1$ (resp.
$\widetilde{\lambda}_0$ and $\widetilde{\lambda}_1$) the extremal leaves of $F_{\widetilde{\lambda}}$ (resp. $F'_{\widetilde{\lambda}}$), i.e. the leaves that satisfy
$F_{\widetilde{\lambda}}=B_{\Lmr}(\widetilde{\ell}_0,\widetilde{\ell}_1)$
and $F_{\widetilde{\lambda}}'=B_{\widetilde{\Lambda}'}(\widetilde{\lambda}_0,\widetilde{\lambda}_1)$.
Then there exists a unique map $\phi_{\widetilde{\lambda}}
:F_{\widetilde{\lambda}}\to F_{\widetilde{\lambda}}'$ such that for every  
$\widetilde{\ell}\in F_{\widetilde{\lambda}}$, we have $\nu'(B_{\widetilde{\Lambda}'}(\phi_{\widetilde{\lambda}}(\widetilde{\ell}_0),
\phi_{\widetilde{\lambda}}(\widetilde{\ell}))=
\nu_{\mutilde_{[q]}}(B_{\Lqr}(\widetilde{\ell}_0,\widetilde{\ell}))$.
We denote by $\phi$ the map from $\Lqr$ to $\widetilde{\Lambda}'$ equal to $\varphi_{|\Lqr}$ outside the preimages of the atoms of $\nu_{\mutilde_m}$ 
and equal to $\phi_{\widetilde{\lambda}}$ on the sets $F_{\widetilde{\lambda}}$ where $\widetilde{\lambda}$ is an atom of $\nu_{\mutilde_m}$.
Then, by construction, for every pair of leaves $\widetilde{\ell}_0$ and $\widetilde{\ell}_1$ of $\Lqr$, we have $B_{\widetilde{\Lambda}'}(\phi(\widetilde{\ell}_0),
\phi(\widetilde{\ell}_1))=
\phi(B_{\Lqr}(\widetilde{\ell}_0,\widetilde{\ell}_1))$. Hence the map $\phi$ defines a map 
$\phi_T:(T_{[q]},d_{T_{[q]}})\to(T_{m},d_{T_{m}})$, where 
$(T_{[q]},d_{T_{[q]}})$ and $(T_{m},d_{T_{m}})$ are respectively the tree associated to 
$(\Lqr,\mutilde_{[q]})$ and the dual tree to $(\Lmr,\mutilde_m)$.

\blemm \label{isometrieequivariante}
The map $\phi_T:(T_{[q]},d_{T_{[q]}})\to(T_{m},d_{T_{m}})$ is a $\grperevet$-equivariant isometry.
\elemm

\dem If $\widetilde{\ell}_0$ and $\widetilde{\ell}_1$ are two leaves of $\Lqr$, we have 
$\phi(B_{\Lqr}(\widetilde{\ell}_0,\widetilde{\ell}_1))=B_{\widetilde{\Lambda}'}(\phi(\widetilde{\ell}_0),\phi(\widetilde{\ell}_1))$
and since
$\nu_{\mutilde_m}=\varphi_*\nu_{\mutilde_{[q]}}$, the map $\phi_T$ is isometric. Moreover, $\phi$ is surjective hence, by taking the quotient, 
so is $\phi_T$ and $\phi_T$ is an isometry. Finally, as $\varphi$ is, the map
$\phi$ is $\grperevet$-equivariant and by taking the quotient, so is $\phi_T$.\cqfd

\bibliographystyle{alphanum}
\bibliography{biblio}{}
Département de mathématique, UMR 8628 CNRS, Université Paris-Sud, Bât. 430, F-91405 Orsay Cedex, FRANCE. 
Bureau : 16.

{\it thomas.morzadec@math.u-psud.fr}
\end{document}